\documentclass{amsart}
\usepackage{amssymb}

\newcommand{\bbN}{{\mathbb{N}}}
\newcommand{\bbR}{{\mathbb{R}}}

\newcommand{\bbC}{{\mathbb{C}}}

\newcommand{\calB}{{\mathcal B}}

\newcommand{\calD}{{\mathcal D}}
\newcommand{\calE}{{\mathcal E}}

\newcommand{\calH}{{\mathcal H}}

\newcommand{\calK}{{\mathcal K}}
\newcommand{\calL}{{\mathcal L}}

\newcommand{\calP}{{\mathcal P}}

\newcommand{\no}{\nonumber}
\newcommand{\lb}{\label}

\newcommand{\wti}{\widetilde  }

\newcommand{\loc}{\text{\rm{loc}}}
\newcommand{\spec}{\text{\rm{spec}}}

\newcommand{\ran}{\text{\rm{ran}}}
\newcommand{\ind}{\text{\rm{ind}}}
\newcommand{\dom}{\text{\rm{dom}}}

\newcommand{\bi}{\bibitem}
\newcommand{\hatt}{\widehat}
\newcommand{\beq}{\begin{equation}}
\newcommand{\eeq}{\end{equation}}
\newcommand{\ba}{\begin{align}}
\newcommand{\ea}{\end{align}}

\renewcommand{\Re}{\text{\rm Re}}
\renewcommand{\Im}{\text{\rm Im}}
\renewcommand{\ln}{\text{\rm ln}}

\DeclareMathOperator{\tr}{tr}

\DeclareMathOperator*{\nlim}{n-lim}
\DeclareMathOperator*{\slim}{s-lim}

\allowdisplaybreaks

\numberwithin{equation}{section}

\newtheorem{theorem}{Theorem}[section]
\newtheorem{lemma}[theorem]{Lemma}

\newtheorem{hypothesis}[theorem]{Hypothesis}
\theoremstyle{definition}
\newtheorem{definition}[theorem]{Definition}

\theoremstyle{remark}
\newtheorem{remark}[theorem]{Remark}

\begin{document}
\title[Applications of the Spectral Shift Operator]
{Some Applications of the Spectral Shift  Operator}
\author[Gesztesy and Makarov]{Fritz Gesztesy and
Konstantin A.~Makarov}
\address{Department of Mathematics,
University of
Missouri, Columbia, MO
65211, USA}
\email{fritz@math.missouri.edu}
\urladdr{http://www.math.missouri.edu/people/fgesztesy.html}
\address{Department of Mathematics, University of
Missouri, Columbia, MO
65211, USA}
\email{makarov@azure.math.missouri.edu}

\subjclass{Primary 47B44, 47A10; Secondary 47A20, 47A40}


\begin{abstract}
The recently introduced concept of a spectral shift operator 
is applied in several instances. Explicit applications 
include Krein's trace formula for pairs of self-adjoint 
operators, the Birman-Solomyak spectral averaging formula 
and its operator-valued extension, and an abstract approach  
to trace formulas based on perturbation theory and the 
theory of self-adjoint extensions of symmetric operators.

\end{abstract}

\maketitle
\section{Introduction}\lb{s1}
The concept of a spectral shift function, historically first
 introduced by  I.~M.~Lifshits \cite{Li52}, \cite{Li56} 
and then developed into a powerful spectral theoretic 
tool by M.~Krein \cite{Kr62}, \cite{Kr83},
\cite{KJ81}, attracted considerable attention in the past
due to its widespread applications in a variety of fields
including scattering theory, relative index theory, spectral
averaging and its application to localization
properties of random Hamiltonians, eigenvalue
counting functions and spectral asymptotics, semi-classical
approximations, and trace formulas for
one-dimensional Schr\"odinger and Jacobi operators. For an  
extensive bibliography in this connection we refer to 
\cite{GMN98}; detailed  reviews on the spectral shift
function and its applications were published by
Birman and Yafaev \cite{BY93}, \cite{BY93a} in 1993.

The principal aim of this paper is to follow up on our recent 
paper \cite{GMN98}, which was devoted to the introduction of a 
spectral shift operator $\Xi(\lambda,H_0,H)$ for a.e. 
$\lambda\in\bbR,$ associated
with a pair of self-adjoint operators $H_0,$ $H=H_0+V$ with
$V\in\calB_1(\calH)$ ($\calH$ a complex separable Hilbert
space). In the special cases of sign-definite perturbations
$V\geq 0$ and $V\leq 0,$ $\Xi(\lambda,H_0,H)$ turns out to
be a trace class operator in $\calH,$ whose trace coincides
with Krein's spectral shift function $\xi(\lambda,H_0,H)$
for the pair $(H_0,H).$ While the special case $V\geq 0$ has
previously been studied by Carey \cite{Ca76}, our aim in 
\cite{GMN98} was to treat the case of general 
interactions $V$
by separately introducing the positive and negative parts
$V_\pm=(|V|\pm V)/2$ of $V.$ In general, if $V$ is not
sign-definite, then
$\Xi(\lambda,H_0,H)$ (naturally associated with \eqref{3.5}) 
is
not necessarily of trace class. However, we introduced trace
class operators
$\Xi_\pm(\lambda)$ corresponding to $V_\pm,$ 
acting in
distinct Hilbert spaces $\calH_{\pm},$ such that
\begin{equation}
\xi(\lambda,H_0,H)=\tr_{\calH_+}(\Xi_+(\lambda))
- \tr_{\calH_-}(\Xi_-(\lambda)) \text{ for a.e. }
\lambda\in\bbR. \lb{1.1}
\end{equation}
(An alternative approach to $\xi(\lambda,H_0,H),$ which
does not rely on separately introducing $V_+$ and $V_-,$
will be discussed elsewhere \cite{GM99}.)

Our main techniques are based on operator-valued Herglotz
functions (continuing some of our recent investigations
in this area \cite{GKMT98}, \cite{GMT98}, \cite{GT97})
and especially, on a detailed study of logarithms of Herglotz
operators in Section~\ref{s2} following the treatment in 
\cite{GMN98}. In
Section~\ref{s3} we  introduce the
spectral shift operator $\Xi(\lambda,H_0,H)$ associated with
the pair $(H_0,H)$ and relate it to Krein's spectral
shift function $\xi(\lambda,H_0,H)$ and his celebrated
trace formula \cite{Kr62}. Finally, Section~\ref{s4} 
provides various applications of this formalism including
spectral averaging originally due to Birman and Solomyak
\cite{BS75},  its operator-valued generalization first 
discussed
in
\cite{GMN98}, connections with the scattering operator, 
and an 
abstract version of a trace formula, combining
perturbation theory and the theory of self-adjoint
extensions.

\section{Logarithms of Operator-Valued Herglotz Functions} 
\lb{s2}
The principal purpose of this section is to recall the 
basic properties of logarithms and
associated representation theorems for operator-valued 
Herglotz functions following the treatment in \cite{GMN98}.

In the following $\calH$ denotes a complex separable Hilbert 
space
with scalar product $(\,\cdot, \, \cdot)_{\calH}$
(linear in the second factor) and norm
$\vert\vert \cdot \vert\vert_\calH,$  $I_{\calH}$ the identity
operator in $\calH$,
$\calB(\calH)$ the Banach space of bounded linear
operators defined on $\calH$, $\calB_p(\calH), \,\, p\ge 1$
the standard Schatten-von Neumann ideals of
 $\calB(\calH)$ (cf., e.g., \cite{GK69}, \cite{Si79}) and
$\bbC_+ $ (respectively, $\bbC_-$) the open complex upper
(respectively, lower) half-plane. Moreover, real and
imaginary parts of a bounded operator $T\in \calB(\calH)$ are
defined as  usual by
$\Re (T)=(T+T^*)/2$,
$\Im (T)=(T-T^*)/(2i).$

\begin{definition} \lb{d2.1}
$M:\bbC_+\to \calB(\calH)$ is called an
 {\it operator-valued  Herglotz function}
 (in short,
a {\it Herglotz operator}) if $M$ is analytic on $\bbC_+$ and
$\Im (M(z))\ge 0$ for all
$z\in \bbC_+$.
\end{definition}
\begin{theorem}\lb{t2.2} \mbox{\rm (Birman and Entina 
\cite{BE67},
de Branges \cite{de62}, Naboko \cite{Na87}--\cite{Na90}.)}
Let \linebreak $M:\bbC_+\to \calB(\calH)$ be a Herglotz 
operator.

\noindent(i) Then there exist bounded self-adjoint operators
 $A=A^*\in \calB(\calH),$
 $0\le B\in \calB(\calH)$, a Hilbert space
$ \calK\supseteq \calH,$ a self-adjoint
 operator $L=L^*$ in
$ \calK$, a bounded nonnegative operator
$0\le R\in \calB(\calK)$ with $R\vert_{\calK
\ominus \calH}=0$ such that
\begin{subequations}\lb{2.1}
\begin{align}
M(z)&=A+Bz+R^{1/2} 
(I_{\calK}+zL)(L-z)^{-1}R^{1/2}\vert_{\calH} 
\label{2.1a} \\
&= A+(B+R\vert_{\calH})z+(1+z^2)R^{1/2}(L-z)^{-1}
R^{1/2}\vert_{\calH}.
 \label{2.1b}
\end{align}
\end{subequations}

\noindent(ii) Let $p\ge1$. Then $M(z)\in \calB_p(\calH)$ 
for all
$z\in \bbC_+$ if and
only if $M(z_0)\in \calB_p(\calH)$ for some $z_0\in \bbC_+.$
In this case necessarily $A,B,R\in \calB_p(\calH).$

\noindent(iii) Let $M(z) \in \calB_1(\calH)$ for some (and
 hence for all) $z\in \bbC_+$. Then $M(z)$
 has normal boundary values $M(\lambda+i0)$ for (Lebesgue)
a.e. $\lambda\in \bbR$ in
every $\calB_p(\calH)$-norm, $p>1$. Moreover, let
$\{E_L(\lambda)\}_{\lambda\in\bbR}$ be the family of 
orthogonal spectral
projections of $L$ in $\calK.$ Then 
$R^{1/2}E_L(\lambda)R^{1/2}$ is
$\calB_1(\calH)$-differentiable for a.e. $\lambda\in\bbR$ and
denoting the derivative by 
$d(R^{1/2}E_L(\lambda)R^{1/2})/d\lambda,$
$\Im(M(z))$ has normal boundary values $\Im(M(\lambda+i0))$ 
for
a.e. $\lambda\in\bbR$ in $\calB_1(\calH)$-norm given by
\begin{equation}
\lim_{\varepsilon\downarrow 0}\|\pi^{-1}\Im(M(\lambda+
i\varepsilon))
- d(R^{1/2}E_L(\lambda)
R^{1/2}|_\calH)/d\lambda\|_{\calB_1(\calH)}
=0 \text{ a.e. } \lb{2.2a}
\end{equation}
\end{theorem}

Originally, the existence of normal limits $M(\lambda+i0)$ for
a.e. $\lambda\in\bbR$ in $\calB_2(\calH)$-norm, in the special
case $A=0,$ $B=-R|_\calH,$ assuming $M(z)\in\calB_1(\calH),$ 
was
proved by de Branges \cite{de62}
in 1962. (The more general case in \eqref{2.1} can easily be
reduced to this special case.) In his paper \cite{de62}, 
de Branges
also proved the existence of normal limits 
$\Im(M(\lambda+i0))$
for a.e. $\lambda\in\bbR$ in $\calB_1(\calH)$-norm and 
obtained
\eqref{2.2a}.
These results and their implications on stationary scattering 
theory
were subsequently studied in detail by Birman and Entina
\cite{BE64}, \cite{BE67}. (Textbook representations of this
material can also be found in \cite{BW83}, Ch.~3.)

Invoking the family of orthogonal spectral projections
$\{E_L(\lambda)\}_{\lambda\in \bbR}$ of $L$, \eqref{2.1}
then yields the  familiar representation
\begin{equation}\lb{2.2}
M(z)= A+Bz+
\int_{\bbR} (1+\lambda^2)
d(R^{1/2}E_L(\lambda)
 R^{1/2}\vert_{\calH})
  ((\lambda-z)^{-1}
-\lambda(1+\lambda^2)^{-1}) ,
\end{equation}
where for our purpose it suffices to interpret
the integral in \eqref{2.2} in the weak sense.
Further results on representations of the type \eqref{2.2} 
can be found in
\cite{Br71}, Sect.~I.4  and  \cite{Sh71}.

Since we are interested in logarithms of Herglotz operators, 
questions of
 their invertibility naturally arise.
The following result clarifies the situation.
\begin{lemma}\lb{l2.3} \mbox{\rm (\cite{GMN98}.)}
Suppose $M$ is a Herglotz operator with values in 
$\calB(\calH)$.
 If $M(z_0)^{-1}\in \calB(\calH)$ for some $z_0\in \bbC_+$ 
then
$M(z)^{-1}\in \calB(\calH)$ for all $z\in \bbC_+$.
\end{lemma}

Concerning boundary values at the real axis we also recall

\begin{lemma} \lb{l2.3b} \mbox{\rm (Naboko \cite{Na93}.)} 
Suppose
$(M-I_\calH)$ is a Herglotz operator with values in 
$\calB_1(\calH).$
Then the boundary values $M(\lambda+i0)$ exist for a.e.
$\lambda\in\bbR$ in $\calB_p(\calH)$-norm, $p>1$ and
$M(\lambda+i0)$ is a Fredholm operator for a.e. 
$\lambda\in\bbR$
with index zero a.e.,
\begin{equation}
\ind(M(\lambda+i0)) =0 \text{ for a.e. } \lambda\in\bbR.
\lb{2.13s}
\end{equation}
Moreover,
\begin{equation}
\ker(M(\lambda+i0)) =\ker(M(i)) = (\ran(M(\lambda+i0)))^\bot
\text{ for a.e. } \lambda \in\bbR. \lb{2.13t}
\end{equation}
In addition, if $M(z_0)^{-1}\in\calB(\calH)$ for some 
{\rm (}and hence
for all{\rm \,)} $z_0\in\bbC_+,$ then
\begin{equation}
M(\lambda+i0)^{-1}\in\calB(\calH) \text{ for a.e. }
\lambda\in\bbR. \lb{2.13u}
\end{equation}
\end{lemma}

Next, let $T$ be a bounded {\it dissipative} operator, 
that is,
\begin{equation}\lb{2.14}
T\in \calB(\calH), \quad \Im (T) \ge 0.
\end{equation}
In order to define the logarithm of $T$
we use the integral representation
\begin{equation}
\log (z) =-i \int_0^\infty d \lambda \,
 ((z+i\lambda)^{-1}- (1+i\lambda)^{-1}), \quad
z \neq-iy, \, y\geq 0, \lb{2.15}
\end{equation}
with a cut along the negative imaginary $z$-axis.
We use the symbol $\log(\cdot)$ in \eqref{2.15} in order to 
distinguish it
from the integral representation
\begin{equation}\lb{2.15a}
\ln(z) = \int_{-\infty}^0 d\lambda \, ((\lambda - z)^{-1} -
\lambda(1+\lambda^2)^{-1}), \quad z\in
\bbC\backslash (-\infty,0]
\end{equation}
with a cut along the negative real axis. It is easily 
verified
that $\log(\cdot)$ and 
$\ln(\cdot)$ coincide
for $z\in\bbC_+.$ In particular,
one verifies that \eqref{2.15} and
\eqref{2.15a} are Herglotz functions,
 that is, they are analytic in $\bbC_+$
 and
\begin{equation}\lb{2.16}
0 < \Im (\log (z)) , \, \, \Im (\ln (z)) \, < \pi, \quad 
z\in \bbC_+.
\end{equation}
\begin{lemma}\lb{l2.4} \mbox{\rm (\cite{GMN98}.)}
Suppose $T\in \calB(\calH)$ is dissipative
and $T^{-1}\in \calB(\calH)$. Define
\begin{equation}\lb{2.17}
\log (T)=-i\int_0^\infty  d \lambda \,
((T+i\lambda)^{-1}-(1 +i\lambda)^{-1}I_\calH)
\end{equation}
in the sense of a $\calB(\calH)$-norm convergent Riemann 
integral.
Then

\noindent (i) $\log (T)\in \calB(\calH).$

\noindent (ii) If
$T=zI_{\calH},$ $z \in \bbC_+,$
 then $\log(T) =\log (z) I_{\calH}.$

\noindent (iii) Suppose
$\{ P_n\}_{n\in \bbN} \subset \calB(\calH)$ is a
family of  orthogonal finite-rank projections
in $\calH$ with $\slim_{n \to \infty}P_n=I_{\calH}.$ Then
$$
\slim_{n \to \infty}((I_{\calH}-P_n)+P_nTP_n)= T
$$
and
\begin{align}
&\slim_{n \to \infty}\log ((I_{\calH}- P_n) + 
P_n (T+i\varepsilon)P_n)
\no \\
&
=
\slim_{n \to \infty}P_n(\log (P_n (T+
i\varepsilon)P_n\vert_{P_n\calH})P_n =
 \log
(T+i\varepsilon I_{H}), \quad \varepsilon >0.\no
\end{align}

\noindent (iv) $
\nlim_{\varepsilon \downarrow 0 }\log (T+
i\varepsilon I_{\calH}) = \log (T).$

\noindent (v) $e^{\log(T)}=T.$
\end{lemma}

\begin{lemma} \mbox{\rm (\cite{GMN98}.)} \lb{l2.5} 
Suppose $T\in \calB(\calH)$ is 
dissipative
and $T^{-1}\in \calB(\calH)$. Let $L$ be the minimal 
self-adjoint dilation
of $T$ in the Hilbert space $\calK\supseteq \calH$.
Then
\begin{equation}\lb{n2.25}
\Im (\log (T))= \pi P_{\calH} E_L(( -\infty,0))\vert_{\calH},
\end{equation}
where $P_{\calH}$ is the orthogonal projection in $\calK$ 
onto $\calH$ and
$\{ E_L(\lambda)\}_{\lambda\in \bbR}$ is the family of 
orthogonal
 spectral projections of $L$ in $\calK$.
In particular,
\begin{equation}\lb{n2.26}
0\le \Im (\log (T))\le \pi I_{\calH}.
\end{equation}
\end{lemma}

Combining Lemmas~\ref{l2.3} and \ref{l2.5}  one obtains the
following result.
\begin{lemma}\lb{l2.6} \mbox{\rm (\cite{GMN98}.)}
Suppose $M:\bbC_+\longrightarrow \calB(\calH)$ is
a Herglotz operator and assume that 
$M(z_0)^{-1}\in \calB(\calH)$ for some
{\rm (}and hence for all{\rm \,)} $z_0\in \bbC_+$. Then
$\log (M):\bbC_+ \longrightarrow \calB(\calH)$ is a Herglotz 
operator
and
\begin{equation}\lb{2.32}
0\le \Im (\log (M(z)))\le \pi I_\calH, \quad z\in \bbC_+.
\end{equation}
\end{lemma}

\begin{theorem}\lb{t2.8} \mbox{\rm (\cite{GMN98}.)}
Suppose $M:\bbC_+ \longrightarrow \calB(\calH)$ is a Herglotz 
operator
and $M(z_0)^{-1} \in \calB(\calH)$ for some {\rm (}and hence
for all{\rm  \,)} $z_0\in\bbC_+.$ Then
there exists a family of bounded self-adjoint
weakly {\rm (}Lebesgue{\rm \,)} measurable operators
 $\{\Xi(\lambda) \}_{\lambda\in \bbR}\subset \calB(\calH),$
\begin{equation}\lb{2.38}
0\le \Xi(\lambda)\le I_\calH \text{ for a.e. }  
\lambda\in \bbR
\end{equation}
such that
\begin{equation}\lb{2.39}
\log(M(z))=C+
\int_\bbR d \lambda \, \Xi(\lambda)
((\lambda-z)^{-1}-\lambda(1+\lambda^2)^{-1}),
 \quad z\in \bbC_+
\end{equation}
the integral taken in the weak sense, where $C=C^* \in 
\calB(\calH).$
 Moreover, if $\Im (\log$ $(M(z_0)))\in \calB_1(\calH)$ for
 some {\rm (}and hence for all{\rm \,)}
 $z_0\in \bbC_+$, then
\begin{align}
&
0\le \Xi(\lambda)\in \calB_1(\calH)
 \text{ for a.e. }  \lambda\in \bbR,
\lb{2.40} \\
&
0\le \tr_{\calH}(\Xi (\cdot))\in 
L_{\rm{loc}}^1(\bbR;d\lambda), \quad
\int_\bbR d\lambda \, (1+\lambda^2)^{-1}
\tr_\calH (\Xi(\lambda))<\infty,
\lb{2.42}
\end{align}
and
\begin{equation}\lb{2.43}
\tr_\calH (\Im (\log (M(z))))=\Im (z) \int_\bbR
d\lambda \,
\tr_{\calH} (\Xi(\lambda))
\vert \lambda-z\vert^{-2} , \quad z\in \bbC_+.
\end{equation}
\end{theorem}

\begin{remark}\lb{r2.9}
For simplicity we only focused on dissipative operators.
Later  we will
also encounter
operators  $S\in \calB(\calH)$ with $-S$ dissipative, 
that is,
$\Im (S) \le 0$ (cf.~\eqref{3.13b}).
In this case $S^*$ is dissipative and one can simply define 
$\log(S)$ by
\begin{equation}
\log(S) = (\log(S^*))^*,   \lb{2.55a}
\end{equation}
with $\log(S^*)$ defined as in \eqref{2.17}.
 Moreover,
\begin{align}
&\log (\hatt M(z))=\hatt C -
\int_\bbR d \lambda \, \hatt \Xi(\lambda)
((\lambda-z)^{-1}-\lambda(1+\lambda^2)^{-1}), \quad
z\in \bbC_+, \lb{2.52} \\
&\widehat C =\widehat C^* \in \calB(\calH)
\text{ and } 0 \le \hatt \Xi(\lambda )\le I_\calH  
\text{ for a.e. }
\lambda\in \bbR, \lb{2.52a}
\end{align}
whenever $\hatt M$ is analytic in $\bbC_+$ and
 $\Im (\hatt M(z))\le 0,$ $z\in \bbC_+.$
\end{remark}
\begin{remark}\lb{r2.10}
Theorem~\ref{t2.8} represents the operator-valued 
generalization of the
exponential Herglotz representation for scalar Herglotz 
functions studied in
detail
by Aronszajn and Donoghue
\cite{AD56} (see also Carey and Pepe
\cite{CP73}). Prior to our proof of Theorem~\ref{t2.8} 
in \cite{GMN98}, Carey \cite{Ca76} 
considered the case $M(z)=I_\calH +K^*(H_0-z)^{-1}K$ 
in 1976 and established
\begin{equation}\lb{2.56}
M(z)=\exp \bigg (\int_\bbR d\lambda \,\Xi (\lambda) 
(\lambda-z)^{-1}
\bigg )
\end{equation}
for  a summable operator function $\Xi (\lambda)$, 
$0\le \Xi (\lambda)\le I_\calH$. 
Carey's proof is different from ours and
does not utilize the integral representation \eqref{2.17} 
for logarithms.
\end{remark}

\section{The Spectral Shift Operator}\lb{s3}

The main purpose of this section is to recall the concept
of a spectral shift operator (cf.~Definition~\ref{d3.4a}) 
as developed in \cite{GMN98}.

Suppose $\calH$ is a complex separable Hilbert space 
and assume 
the following
hypothesis for the remainder of this section.
\begin{hypothesis}\lb{h3.1}
Let  $H_0$ be a self-adjoint operator in $\calH$
with domain $\dom (H_0)$, $J$ a bounded self-adjoint 
operator
with $J^2=I_{\calH}$, and $K\in \calB_2(\calH)$ a 
Hilbert-Schmidt operator.
\end{hypothesis}

Introducing
\begin{equation}\lb{3.1}
V=KJK^*
\end{equation}
we define the self-adjoint operator
\begin{equation}\lb{3.2}
H=H_0+V, \quad \dom(H)=\dom(H_0)
\end{equation}
in $\calH$.

Given Hypothesis~\ref{h3.1} we decompose $\calH$ and $J$ 
according to
\begin{equation}\lb{3.3}
J=\begin{pmatrix} I_+ & 0\\ 0&
-I_-\end{pmatrix}, \quad \calH=\calH_+\oplus \calH_-,
\end{equation}
\begin{equation}\lb{3.4}
J_+=\begin{pmatrix} I_+ & 0\\ 0&
0\end{pmatrix}, \quad
J_-=\begin{pmatrix} 0& 0\\ 0&
I_-\end{pmatrix}, \quad
J=J_+-J_-,
\end{equation}
with $I_\pm$ the identity operator in $\calH_\pm$. Moreover, 
we introduce
the following bounded operators
\begin{align}
\Phi (z)
&=J+K^*(H_0-z)^{-1}K:\calH \rightarrow \calH, \lb{3.5}\\
\Phi_+(z)
&=I_++J_+ K^*(H_0-z)^{-1}K\vert_{\calH_+}:
\calH_+ \rightarrow\calH_+, \lb{3.6} \\
\widetilde \Phi_-(z)
&=I_- - J_-K^*(H_+-z)^{-1}K\vert_{\calH_-}:
\calH_- \rightarrow \calH_-, \lb{3.7}
\end{align}
for $z\in \bbC\backslash \bbR,$ where
 \begin{equation}\lb{3.7'}
V_+=KJ_+K^*,
\end{equation}
\begin{equation}\lb{3.7''}
H_+=H_0+V_+, \quad \dom (H_+)=\dom (H_0).
\end{equation}
\begin{lemma}\lb{l3.2} \mbox{\rm (\cite{GMN98}.)}
Assume Hypothesis~\ref{h3.1}. Then $\Phi$, $\Phi_+,$
 and $-\widetilde \Phi_-$ are Herglotz operators in
$\calH,$  $\calH_+,$ and
$\calH_-,$ respectively. In addition {\rm (}$z\in
\bbC\backslash \bbR${\rm )},
\begin{align}
\Phi (z)^{-1}
&=J-JK^*(H-z)^{-1}KJ, \lb{3.8} \\
\Phi_+ (z)^{-1}
&=I_+-J_+K^*(H_+-z)^{-1}K\vert_{\calH_+}, \lb{3.9} \\
\widetilde \Phi_- (z)^{-1}
&=I_- +J_-K^*(H-z)^{-1}K\vert_{\calH_-}.
\lb{3.10}
\end{align}
\end{lemma}

We also recall the following result (cf.~\cite{GMN98}).

\begin{lemma}\lb{l3.3}
Assume Hypothesis~\ref{h3.1} and $\bbC\backslash\bbR$. 
Then
\begin{subequations} \lb{3.13}
\begin{align}
\tr_\calH((H_0-z)^{-1}-(H_+-z)^{-1})
&=d \tr_{\calH_+} ( \log (\Phi_+(z)))/dz, \lb{3.13a} \\
\tr_\calH((H_+-z)^{-1}-(H-z)^{-1})
&=d \tr_{\calH_-} ( \log (\widetilde\Phi_-(z)))/dz. 
\lb{3.13b}
\end{align}
\end{subequations}
\end{lemma}

Next, applying Theorem~\ref{t2.8} and Remark~\ref{r2.9} to
$\Phi_+(z)$ and $\widetilde \Phi_-(z)$ one infers the 
existence of two families
of bounded operators
$\{\Xi_\pm (\lambda)\}_{\lambda \in \bbR}$ defined for 
(Lebesgue) a.e.
$\lambda\in \bbR$ and
 satisfying
\begin{align}
&0\le \Xi_{\pm}(\lambda)\le I_\pm, \quad
\Xi_\pm (\lambda) \in \calB_1(\calH_\pm)
\text{ for a.e. } \lambda\in \bbR, \lb{3.16} \\
&
\vert\vert
\Xi_\pm(\cdot)\vert\vert_1\in L^1(\bbR;(1+
\lambda^2)^{-1} d\lambda)
\no
\end{align}
and
\begin{subequations} \lb{3.17}
\begin{align}
\log (\Phi_+(z))
&=
\log (I_++J_+K^*(H_0-z)^{-1}K\vert_{\calH_+})
\no\\
&=
C_+ +
\int_\bbR d\lambda \, \Xi_+(\lambda)
((\lambda-z)^{-1}-\lambda (1+\lambda^2)^{-1}),
\lb{3.17a} \\
\log (\widetilde \Phi_-(z))
&=
\log (I_--J_-K^*(H_+-z)^{-1}K\vert_{\calH_-})
\no \\
&=C_- -
\int_\bbR d\lambda  \,\Xi_-(\lambda)
((\lambda-z)^{-1}-\lambda (1+\lambda^2)^{-1})
\lb{3.17b}
\end{align}
\end{subequations}
for $z\in \bbC\backslash \bbR, $ with $C_\pm =C_\pm^*\in
\calB_1(\calH)$.

Equations \eqref{3.17} motivate the following
\begin{definition} \lb{d3.4a}
$\Xi_+(\lambda)$ (respectively, $\Xi_-(\lambda)$) is 
called the
{\it spectral shift operator} associated with
$\Phi_+(z)$ (respectively, $\widetilde \Phi_-(z)$). 
Alternatively, 
we will
refer to
$\Xi_+(\lambda)$ as the spectral shift operator
associated with the pair $(H_0,H_+)$ and occasionally
use the notation $\Xi_+(\lambda,H_0,H_+)$ to stress the 
dependence on
$(H_0,H_+),$ etc.
\end{definition}

Moreover, we introduce
\begin{equation}\lb{3.19}
\xi_\pm (\lambda) =\tr_{\calH_\pm}( \Xi_{\pm}(\lambda) ), 
\quad
0\le \xi_\pm\in L^1(\bbR;(1+\lambda^2)^{-1}d\lambda)
\text{ for a.e. } \lambda \in \bbR.
\end{equation}

Actually, taking into account the simple  behavior of
$\Phi_+(iy)$ and
$\widetilde\Phi_-(iy)$ as $
\vert y\vert \to \infty,$ one can improve
\eqref{3.17a} and \eqref{3.17b} as follows.
\begin{lemma}\lb{l3.4} \mbox{\rm (\cite{GMN98}.)}
Assume Hypothesis~\ref{3.1}
 and define $\xi_\pm$ as in \eqref{3.19}. Then
\begin{equation}\lb{3.19''}
0\le \xi_\pm \in L^{1}(\bbR; d\lambda),
\end{equation}
and \eqref{3.17a} and \eqref{3.17b} simplify to
\begin{subequations} \lb{3.25}
\begin{align}
\log(\Phi_+(z))
&=\int_\bbR d\lambda \, \Xi_+(\lambda)(\lambda -z)^{-1}, 
\lb{3.25a}
\\
\log(\widetilde\Phi_-(z))
&=-\int_\bbR d\lambda \, \Xi_-(\lambda)(\lambda -z)^{-1}.
\lb{3.25b}
\end{align}
\end{subequations}
Moreover, for a.e. $\lambda\in\bbR,$
\begin{subequations} \lb{3.26a}
\begin{align}
\lim_{\varepsilon\downarrow 0}\|\Xi_+(\lambda)-
\pi^{-1}\Im(\log(\Phi_+(\lambda+i\varepsilon)))\|
_{\calB_1(\calH_+)} &=0, \lb{3.26aa} \\
\lim_{\varepsilon\downarrow 0}\|\Xi_-(\lambda) +
\pi^{-1}\Im(\log(\wti \Phi_-(\lambda+i\varepsilon)))\|
_{\calB_1(\calH_-)} &=0. \lb{3.26ab}
\end{align}
\end{subequations}
\end{lemma}
\begin{proof}
For convenience of the reader we reproduce here the proof 
first presented in \cite{GMN98}. It suffices to consider
$\xi_+(\lambda)$ and
$\Phi_+(z).$  Since
\begin{equation}\lb{3.19b}
\vert\vert
\log (\Phi_+(y) )\vert\vert_1=O
(\vert y
\vert^{-1}) \text{ as }
\vert y\vert \to  \infty
\end{equation}
by the Hilbert-Schmidt hypothesis on $K$ and the fact
$\vert\vert (H_0-iy)^{-1}\vert\vert =O
(\vert y
\vert^{-1}) $
as $\vert y\vert \to \infty$, the scalar Herglotz function
$\tr_{\calH_+}(\log (\Phi_+(z)))$ satisfies
\begin{equation}\lb{3.19c}
\vert\tr_{\calH_+}(\log (\Phi_+(z)))\vert
 =O
(\vert y
\vert^{-1})
\text{ as } \vert y\vert \to \infty.
\end{equation}
By standard results (see, e.g., \cite{AD56}, \cite{KK74}),
\eqref{3.19c} yields
\begin{equation}\lb{3.19d}
\tr_{\calH_+}(\log (\Phi_+(z)))=
\int_\bbR d \omega_+(\lambda) (\lambda-z)^{-1}, \quad
z\in \bbC\backslash\bbR,
\end{equation}
where $\omega_+$ is a finite measure,
\begin{equation}\lb{3.19e}
\int_\bbR d \omega_+(\lambda)=-i\lim_{y\uparrow\infty}
(y\tr_{\calH_+}(\log (\Phi_+(z))))<\infty.
\end{equation}
Moreover, the fact that $\Im(\log(\Phi_+(z)))$ is uniformly 
bounded
with respect to $z\in\bbC_+$ yields that $\omega_+$ is 
purely absolutely
continuous,
\begin{equation}\lb{3.19f}
d \omega_+(\lambda)=\xi_+(\lambda)d\lambda,
\quad \xi_+\in L^1(\bbR ; d\lambda),
\end{equation}
where
\begin{align}
\xi_+(\lambda)
&=\pi^{-1} \lim_{\varepsilon \downarrow 0}
(\Im (\tr_{\calH_+}(\log(\Phi_+(\lambda+i\varepsilon)))))
=\tr_{\calH_+}(\Xi_+(\lambda))
\no \\
&
=\pi^{-1}\lim_{\varepsilon \downarrow 0}
(\Im (\log (\text{det}_{\calH_+} (\Phi_+(\lambda+
i\varepsilon)))))
\text{ for a.e. }
\lambda\in \bbR. \lb{3.19g}
\end{align}
In order to prove \eqref{3.26aa} we first observe that
$\Im(\log(\Phi_+(\lambda+i\varepsilon)))$ takes on boundary 
values
$\Im(\log(\Phi_+(\lambda+i0)))$ for a.e. $\lambda\in\bbR$ in
$\calB_1(\calH_+)$-norm by \eqref{2.2a}. Next, choosing an 
orthonormal
system $\{e_n\}_{n\in\bbN}\subset\calH_+,$ we recall that
the quadratic form $(e_n,\Im(\log(\Phi(\lambda+
i0)))e_n)_{\calH_+}$
exists for all $\lambda\in\bbR\backslash\calE_n,$ where 
$\calE_n$
has Lebesgue measure zero. Thus one observes,
\begin{align}
&\lim_{\varepsilon\downarrow 0} (e_m,\Im(\log(\Phi_+(\lambda+
i\varepsilon)))e_n)_{\calH_+}
= (e_m,\Im(\log(\Phi_+(\lambda+i0)))e_n)_{\calH_+} \no \\
&=\pi (e_m,\Xi_+(\lambda)e_n)_{\calH_+} \text{ for }
\lambda\in\bbR\backslash\{\calE_m\cup\calE_n\}. \lb{3.33a}
\end{align}
Let $\calE=\cup_{n\in\bbN} \calE_n,$ then $|\calE|=0$ 
($|\cdot|$
denoting the Lebesgue measure on $\bbR$) and hence
\begin{align}
&(f,\Im(\log(\Phi_+(\lambda+i0)))g)_{\calH_{+}} =\pi(f,
\Xi_+(\lambda)g)_{\calH_+} \lb{3.33b} \\
& \text{for } \lambda\in\bbR\backslash\calE
\text{ and } f,g \in \calD=\text{lin.span}\,
\{e_n\in\calH_+\,|\,n\in\bbN\}.
\no
\end{align}
Since $\calD$ is dense in $\calH_+$ and $\Xi_+(\lambda)\in
\calB(\calH_+)$ one infers
 $\Im(\log(\Phi_+(\lambda+i0)))=\pi\Xi_+(\lambda)$ for a.e.
$\lambda \in \bbR,$ completing the proof.
\end{proof}

Assuming Hypothesis~\ref{h3.1} we define
\begin{equation}\lb{3.20}
\xi(\lambda)=
\xi_+(\lambda)-
\xi_-(\lambda)
 \text{ for a.e. }
\lambda\in \bbR
\end{equation}
and call $\xi(\lambda)$ (respectively, $\xi_+(\lambda),$ 
$\xi_-(\lambda)$) the
spectral shift function associated with the pair $(H_0,H)$ 
(respectively,
$(H_0,H_+),$ $(H_+,H)$), sometimes also denoted by 
$\xi(\lambda,H_0,H),$
etc., to underscore the dependence on the pair involved.

M.~Krein's basic trace formula
\cite{Kr62} is now obtained as follows.
\begin{theorem}\lb{t3.4}
Assume Hypothesis~\ref{h3.1}. Then {\rm (}$z \in \bbC 
\backslash
\{\spec (H_0) \cup \spec(H) \}${\rm )}
\begin{equation}
\tr_\calH((H-z)^{-1}-(H_0-z)^{-1})
=-\int_\bbR d\lambda \, \xi(\lambda) (\lambda-z)^{-2}.
\lb{3.21}
\end{equation}
\end{theorem}
\begin{proof}
Let $z\in \bbC\backslash \bbR$.
By
\eqref{3.19d} and \eqref{3.19f} we infer
\begin{align}
&\tr_{\calH_+}
(\log (\Phi_+(z)))=
\int_\bbR d\lambda \,\xi_+(\lambda)  (\lambda-z)^{-1} ,
\lb{3.22} \\
&
\tr_{\calH_-}
(\log (\widetilde \Phi_-(z)))=
-\int_\bbR d\lambda \, \xi_-(\lambda) (\lambda-z)^{-1} . 
\lb{3.23}
\end{align}
Adding \eqref{3.13a} and \eqref{3.13b}, differentiating
\eqref{3.22} and \eqref{3.23} with respect to $z$
 proves \eqref {3.21} for $z\in \bbC\backslash \bbR$. The 
result extends to all
 $z\in \bbC \backslash \{\spec (H_0) \cup \spec(H) \}$
 by  continuity
of
$((H-z)^{-1}-(H_0-z)^{-1})$ in $\calB_1(\calH)$-norm.
\end{proof}

In particular, $\xi(\lambda)$ introduced in 
\eqref{3.20} is
Krein's original spectral shift function. As
noted in Section~\ref{s2}, the spectral shift operator
$\Xi_+(\lambda)$ in the particular case $V=V_+,$ and its
relation to Krein's spectral shift function 
$\xi_+(\lambda),$
was first studied by Carey \cite{Ca76} in 1976.
\begin{remark}\lb{r3.5}
(i) As shown originally by M.~Krein \cite{Kr62}, the 
trace formula
\eqref{3.21} extends to
\begin{equation}\lb{3.24}
\tr(f(H)-f(H_0))
=\int_{\bbR} d\lambda \,\xi(\lambda) f'(\lambda)
\end{equation}
for appropriate functions $f$. This fact has been
studied by numerous authors and we refer, for instance, to
\cite{BW83}, Ch.~19, \cite{BP98}--\cite{BY93},
\cite{Kr83},
\cite{Kr89}, \cite{Pe85},
\cite{Si75},
\cite{SM94},
\cite{Vo87},
\cite{Ya92}, Ch.~8 and the references therein.

\noindent (ii) Concerning scattering theory for the pair 
$(H_0,H),$ we remark that
$\xi(\lambda)$,  for a.e. $\lambda \in \spec_{ac}(H_0)$
(the absolutely continuous spectrum of $H_0$), is related 
to the scattering
operator at fixed energy $\lambda$ by the Birman-Krein 
formula \cite{BK62},
\begin{equation}\lb{3.25d}
\text{det}_{\calH_\lambda}(S(\lambda, H_0, H))=
e^{-2\pi i \xi(\lambda)}
\text{ for a.e. }
\lambda \in \spec_{ac}(H_0).
\end{equation}
Here
 $S(\lambda, H_0, H)$ denote the fibers in the direct 
integral representation
of the scattering operator
$$
S( H_0, H)=\int_{\spec_{ac}(H_0)}^{\oplus}d\lambda  \,
 S(\lambda, H_0, H)
\text{ in } \calH =
\int_{\spec_{ac}(H_0)}^{\oplus}d\lambda \, 
\calH_\lambda
$$
with respect to the absolutely continuous part 
$H_{0, ac}$ 
of $H_0$.
This fundamental connection, originally due to Birman 
and Krein
\cite{BK62},
is further discussed in
\cite{BW83}, Ch.~19,
\cite{BY93},
\cite{BY93a},
\cite{Ca76},
\cite{Ka78},
\cite{Kr83},
\cite{So93},
\cite{Ya92}, Ch.~8 and the literature cited therein. We 
briefly return to this topic in Lemma~\ref{l4.7}.

\noindent (iii) The standard identity (\cite{GK69}, 
Sect.~IV.3)
\begin{equation}\lb{3.26}
\tr_\calH ((H-z)^{-1}-(H_0-z)^{-1})=
-d \log (\text{det}_\calH (I_\calH+V(H_0-z)^{-1}))/dz
\end{equation}
together with the trace formula \eqref{3.21} yields
the well-known connection between perturbation
 determinants and $\xi(\lambda)$,
 also due to M.~Krein
\cite{Kr62}
\begin{equation}\lb{3.27}
\log (\text{det}_\calH (I_\calH+V(H_0-z)^{-1}))=
\int_\bbR d\lambda \, \xi(\lambda)
 (\lambda -z)^{-1} ,
\end{equation}
\begin{equation}\lb{3.28}
\xi(\lambda)=\lim_{\varepsilon\downarrow 0}\pi^{-1}
\Im (
\log (\text{det}_\calH (I_\calH+V(H_0-(\lambda+ 
i0))^{-1})))
\text{ for a.e. }
\lambda\in \bbR,
\end{equation}
\begin{equation}\lb{3.29}
\tr_\calH(V)=
\int_\bbR d \lambda  \, \xi(\lambda), \quad
\quad
\int_\bbR d \lambda \, \vert \xi(\lambda)\vert \le 
\vert\vert V\vert\vert_1.
\end{equation}
This is discussed in more detail,
for instance, in
\cite{BW83}, Ch.~19,
\cite{BY93},
\cite{Ca76},
\cite{Kr60},
\cite{Kr83},
\cite{KJ81},
\cite{Si75}--\cite{So93},
\cite{Ya92}. Relation \eqref{3.28} and the analog of
\eqref{2.2a} for $d(AE_H(\lambda)B)/d\lambda,$ where 
$H=H_0+V,$
$V=B^*A,$ $A,B\in\calB_2(\calH),$ $V=V^*,$ leads to the 
expression
$$
\xi(\lambda)=(-2\pi i)^{-1}\tr_\calH (
\log (I_\calH-2\pi i(I_\calH
-A(H-\lambda-i0)^{-1}B^*)(d(AE_H(\lambda)B^*)/d\lambda)))
$$
for a.e. $\lambda\in\bbR$ (cf., e.g., \cite{BW83}, 
Sects.~3.4.4 and 19.1.4).
\end{remark}

For numerous additional properties and applications of 
Krein's 
spectral shift function we refer to cite{EP97}, 
\cite{GKS95},  
\cite{JK78}, \cite{Ka78}, \cite{KS98}, \cite{Pu97}, 
\cite{So93} 
and the references  therein.

\section{Some Applications}\lb{s4}

In this section we consider some applications of the 
formalism 
developed in  Sections~\ref{s2} and \ref{s3}.

We start by recalling the simple proof of spectral averaging 
and its relation to Krein's spectral  shift function in
\cite{GMN98}, a circle of ideas originating with Birman and
Solomyak \cite{BS75}.

For this purpose we assume the following.

\begin{hypothesis}\lb{h4.1}
Let $H_0$ be a self-adjoint operator in
$\calH$ with $\dom (H_0)$, and assume
$\{V(s)\}_{s\in \Omega}\subset \calB_1(\calH)$ to be a 
family of 
self-adjoint trace class operators in $\calH$, where 
$\Omega \subseteq
 \bbR$ denotes an open interval. Moreover, suppose that 
$V(s)$
 is continuously differentiable with respect to $s\in\Omega$
in trace norm.
\end{hypothesis}
To begin our discussions we temporarily assume that 
$V(s)\geq 0$, that 
is, we suppose 
\begin{equation} \lb{4.2}
V(s)=K(s)K(s)^*, \quad s\in \Omega
\end{equation}
for some $K(s)\in\calB_2(\calH),$ $s\in\Omega.$ Given 
Hypothesis~\ref{h4.1} we define the self-adjoint operator 
$H(s)$ in $\calH$ by
\begin{equation}
H(s)=H_0+V(s), \quad \dom (H(s))=\dom (H_0),
\quad s\in \Omega. \lb{4.3}
\end{equation}

In  analogy to \eqref{3.5}  and \eqref{3.6} we introduce
in $\calH$ $(s\in \Omega,$ $  z\in \bbC \backslash \bbR$),
\begin{equation}\lb{4.7}
\Phi(z,s)
=I_{\calH}+ K(s)^*(H_0-z)^{-1}K(s)
\end{equation}
and hence infer from Lemma~\ref{l3.2} that
\begin{equation}\lb{4.9}
\Phi (z,s)^{-1}
=I_\calH-K(s)^*(H(s)-z)^{-1}K(s).
\end{equation}
The following is an elementary but useful result needed 
in the context of Theorem~\ref{t4.3}.
\begin{lemma}\lb{l4.2}
Assume Hypothesis~\ref{h4.1} and \eqref{4.2}. Then
$(s\in \Omega, $ $ z \in \bbC\backslash \bbR),$
\begin{equation}\lb{4.11}
d \tr_{\calH}
(\log (\Phi(z,s)))/ds=
\tr_{\calH}(V'(s)(H(s)-z)^{-1}).
\end{equation}
\end{lemma}

Next, applying Lemma~\ref{l3.4} to
$\Phi(z,s)$ in \eqref{4.7} one infers $(s\in \Omega),$
\begin{align}
&\log (\Phi(z,s))=\int d\lambda \, \Xi(\lambda,s)
(\lambda-z)^{-1},
\lb{4.19a} \\
&0\le \Xi (\lambda,s) \le I_\calH, \quad
\Xi (\lambda ,s)\in \calB_1(\calH)
\text{ for a.e. } \lambda \in \bbR, \lb{4.21} \\
&\vert\vert \Xi (\cdot,s)\vert\vert_1 \in L^1(\bbR;d\lambda), 
\no
\end{align}
where $\Xi(\lambda,s)$ is associated with the pair 
$(H_0,H(s)),$ assuming $H(s)\geq H_0,$ $s\in\Omega.$

The principal result on averaging the spectral measure of
$\{E_{H(s)}(\lambda)\}_{\lambda\in \bbR}$ of $H(s)$ 
as proven in \cite{GMN98} then reads as follows.

\begin{theorem}\lb{t4.3} \mbox{\rm (\cite{GMN98}.)}
Assume Hypothesis~\ref{h4.1} and 
$[s_1,s_2]\subset \Omega$. Let $\xi(\lambda,s)$ be the 
spectral shift function associated with the pair 
$(H_0,H(s))$ {\rm (}cf. \eqref{3.20}{\rm )}, where 
$H(s)$ is defined 
by \eqref{4.3} {\rm (}and we no longer suppose 
$H(s)\geq H_0${\rm )}. Then
\begin{equation} \lb{4.25}
\int_{s1}^{s_2} ds \, ( d (\tr_{\calH} 
(V'(s)E_{H(s)}(\lambda) )))=
(\xi(\lambda, s_2)-\xi(\lambda, s_1))d\lambda.
\end{equation}
\end{theorem}

\begin{remark}\lb{r4.4}
 (i) In the special case of averaging over the boundary 
condition
 parameter for half-line Sturm-Liouville operators 
(effectively a
rank-one resolvent perturbation problem), Theorem~\ref{t4.3} 
has
 first been derived by Javrjan
\cite{Ja66},
\cite{Ja71}. The case of rank-one perturbations was 
recently treated
in detail by Simon \cite{Si95}.
The general case of trace class perturbations is due to 
Birman and
 Solomyak
\cite{BS75} using an approach of Stieltjes' double operator 
integrals.
Birman and Solomyak treat the case $V(s)=sV,$ 
$V\in\calB_1(\calH),$ $s\in [0,1].$ A short proof of
\eqref{4.25} (assuming $V'(s)\ge 0$) has recently been 
given by Simon \cite{Si98}.

\noindent (ii) We note that variants  of \eqref{4.25} in  
the context of
one-dimensional
 Sturm-Liouville operators (i.e.,
variants of Javrjan's results in
\cite{Ja66},
\cite{Ja71})
have been repeatedly rediscovered by several authors. In 
particular,
the absolute continuity of averaged spectral measures (with 
respect to
 boundary condition parameters or coupling constants of 
rank-one perturbations)
has been used to prove localization properties
 of one-dimensional random Schr\"odinger operators 
(see, e.g.,
\cite{BS98}, \cite{Ca83},
\cite{Ca84}, \cite{CL90}, Ch.~VIII, \cite{CH94},
\cite{CHM96}, \cite{KS98}--\cite{KS87},
\cite{PF92}, Ch.~V, \cite{Si85}, \cite{Si95}).

\noindent (iii) We emphasize that Theorem~\ref{t4.3} 
applies 
to unbounded
operators (and
hence to random Schr\"odinger operators bounded from below) 
as long as
appropriate relative trace class conditions (either with 
respect to resolvent
or semigroup perturbations) are satisfied.

\noindent (iv) In the special case
$V'(s)\geq 0,$ the  measure 
$$
d(\tr_\calH (V'(s)E_{H(s)}(\lambda))) =
d(\tr_\calH (V'(s)^{1/2}E_{H(s)}(\lambda)V'(s)^{1/2}))
$$
 in \eqref{4.25}
represents a positive measure.
\end{remark}

In the special case of a sign-definite perturbation of $H_0$ 
of the form
$sKK^*,$ one can in fact prove an
operator-valued averaging formula as follows.
\begin{theorem} \lb{t4.5} \mbox{\rm (\cite{GMN98}.)}
Assume Hypothesis~\ref{h3.1} and $J=I_\calH.$ Then
\begin{equation}
\int_0^1 ds \, d(K^*E_{H_0 + sKK^*}(\lambda)K)=
\Xi(\lambda)d\lambda, \lb{4.33}
\end{equation}
where $\Xi(\cdot)$ is the spectral shift operator 
associated with
\begin{equation}
\Phi(z)=I_\calH + K^*(H_0 - z)^{-1}K, 
\quad z\in\bbC\backslash\bbR,
\lb{4.34}
\end{equation}
that is,
\begin{align}
&\log(\Phi(z))= \int_\bbR d\lambda \, \Xi(\lambda) 
(\lambda - z)^{-1},
\quad
z\in\bbC\backslash\bbR, \lb{4.35} \\
& 0\leq \Xi(\lambda) \in \calB_1 (\calH) 
\text{ for a.e. }\lambda\in\bbR,
\quad \|\Xi(\cdot)\|_1\in L^1(\bbR;d\lambda). \lb{4.36}
\end{align}
\end{theorem}
\begin{proof}
An explicit computation shows
\begin{align}
&(\lambda - \Phi(z))^{-1} = -(1-\lambda)^{-1} \no \\
&\times
(I_\calH
-(1-\lambda)^{-1}K^*(H_0+(1-\lambda)^{-1}KK^* - z)^{-1}K) 
\in\calB(\calH)
 \lb{4.37}
\end{align}
for all $\lambda <0.$ Since $\log(\Phi(z))=\ln(\Phi(z))$ for
$z\in\bbC_+$ as a result of analytic continuation, one 
obtains
\begin{align}
\log(\Phi(z)) &= \int_\bbR d\lambda \, \Xi(\lambda) 
(\lambda -z)^{-1}
\no \\
=\ln(\Phi(z)) &= \int_{-\infty}^0 d\lambda \, 
((\lambda -\Phi(z))^{-1} -
\lambda(1+\lambda^2)^{-1}I_\calH) \no \\
&=\int_{-\infty}^0
d\lambda \, (1-\lambda)^{-2} K^*
(H_0+(1-\lambda)^{-1}KK^*-z)^{-1}K \no \\
&= \int_0^1 ds \int_\bbR
d(K^*E_{H_0+sKK^*}(\lambda)K)(\lambda-z)^{-1} \no \\
&=\int_\bbR (\lambda-z)^{-1}\int_0^1
ds \, d(K^*E_{H_0+sKK^*}(\lambda)K) \lb{4.38}
\end{align}
proving \eqref{4.33}. (Here the interchange of the 
$\lambda$ and $s$
integrals follows from Fubini's theorem considering 
\eqref{4.38}
in the weak sense.)
\end{proof}

As a consequence of Theorem~\ref{t4.5} one obtains
\begin{equation}
\int_{s_1}^{s_2} ds\, d(K^* E_{H_0+sKK^*}(\lambda)K) =
\Xi(\lambda,s_2) - \Xi(\lambda,s_1), \lb{4.39}
\end{equation}
where $\Xi(\lambda,s)$ is the spectral shift operator 
associated with
$\Phi(z,s)=I_\calH + sK^*(H_0-z)^{-1}K,$ $s\in [s_1,s_2].$ 
This yields an alternative proof of the following result 
of Carey \cite{Ca76}.

\begin{lemma} \lb{l4.6} \mbox{\rm (\cite{Ca76}.)}
Assume Hypothesis~\ref{h3.1} and $J=I_\calH.$ Then
\begin{equation}\lb{4.40}
K^*K=
\int_{\bbR} d\lambda \,\Xi(\lambda),
\end{equation}
with $\Xi(\lambda)$ given by \eqref{4.35}.
\end{lemma}
\begin{proof}
Given $f\in \calH$ one infers
\begin{align}
\bigg (\int_{\bbR}\int_0^1 ds \, d(K^*E_{H_0 + sKK^*}
(\lambda)K)f,f\bigg )&=
\int_0^1 ds \int_{\bbR}d(K^*E_{H_0 + sKK^*}(\lambda)Kf,f)
\no \\
&=\int_0^1 ds ||Kf||^2 =||Kf||^2 \lb{4.41}
\end{align}
by Fubini's theorem. Combining \eqref{4.33} and \eqref{4.41} 
one concludes
\begin{equation}\lb{4.42}
||Kf||^2=\bigg (\int_{\bbR} d\lambda \,\Xi(\lambda) 
f,f\bigg ).
\end{equation}
Since $K^*K$ and $\int_{\bbR} d\lambda \,\Xi(\lambda)$ are 
both bounded operators and by \eqref{4.42} their quadratic 
forms coincide, one obtains \eqref{4.40}.
\end{proof}

For nonnegative perturbations $V\ge 0$, 
$V\in \calB_1(\calH),$ 
the particular choice $K=K^*=V^{1/2}$ then reconstructs 
$V$ via
\begin{equation} \lb{4.43}
V=\int_{\bbR} d\lambda \,\Xi(\lambda)
\end{equation}
as a consequence of Lemma~\ref{l4.6}.

The case of nonpositive perturbations can be treated in an 
analogous way. The general case of sign indefinite 
perturbations will be considered elsewhere \cite{GM99}.

Next, combining some results of Kato~\cite{Ka78} and 
Carey~\cite{Ca76}, we briefly turn to scattering theory for
the  pair $(H_0,H),$ $H=H_0+KK^*,$ assuming
Hypothesis~\ref{h3.1}  and $J=I_\calH.$ We pick an interval
$(a,b)\subseteq\bbR$  and assume that $H_0$ is spectrally
absolutely continuous on 
$(a,b),$ that is, $E_{H_0}(\lambda)=
E_{H_0,\text{a.c.}}(\lambda),$ $\lambda\in (a,b),$ where 
$\{E_{H_0}(\lambda)\}_{\lambda\in\bbR}$ denotes the family 
of orthogonal spectral projections of $H_0.$ As discussed
in connection with Theorem~\ref{t2.2}, 
\begin{equation}
\Phi(\lambda\pm i0)=\lim_{\varepsilon\downarrow 0} 
\Phi(\lambda\pm i\varepsilon), \lb{4.43a} 
\end{equation}
exist for a.e.~$\lambda\in\bbR$ in $\calB_2(\calH)$-norm
(actually in $\calB_p(\calH)$-norm for all $p>1$)
and
\begin{align}
&A(\lambda)=d(K^*E_{H_0}(\lambda)K)/d\lambda\in
\calB_1(\calH), \lb{4.43b} \\
&\lim_{\varepsilon\downarrow 0}\pi^{-1}\Im(\Phi 
(\lambda+i\varepsilon)) =A(\lambda) \text{ for a.e. } 
\lambda\in\bbR \lb{4.43c}
\end{align}
in $\calB_1(\calH)$-norm (cf., e.g., \cite{BW83},
Sect.~3.4.4), with $\Phi(z)$ defined as in \eqref{4.34}.
Following Kato \cite{Ka78}, one considers $\ran(K)=K\calH$
and defines the semi-inner product 
$(\cdot,\cdot)_{\ran(K),\lambda},$ $\lambda\in (a,b),$ by
\begin{equation}
(x,y)_{\ran(K),\lambda}=(f,A(\lambda)g)_\calH \text{ for
a.e. } \lambda\in (a,b), \, x=Kf, \, y=Kg, \, f,g \in\calH.
\end{equation} \lb{4.43d}
Denoting by $\calK(\lambda)=\overline{(\ran(K);
(\cdot,\cdot)_{\ran(K),\lambda})},$ the completion of
$\ran(K)$ with respect to
$(\cdot,\cdot)_{\ran(K),\lambda},$ the fibers
$S(\lambda,H_0,H),$ $\lambda\in (a,b)$ in the direct
integral representation of the (local) scattering operator 
$S(H_0,H)P_{H_0}((a,b))$ ($P_{H_0}((a,b))$ the
corresponding orthogonal spectral projection of $H_0$
associated with the interval $(a,b)$) then can be identified
with the unitary operator 
\begin{equation}
(I_{\calK(\lambda)}+KK^*(H_0-\lambda-i0)^{-1})^{-1}
(I_{\calK(\lambda)}+KK^*(H_0-\lambda+i0)^{-1}) 
\text{ for a.e. } \lambda\in (a,b) \lb{4.43e} 
\end{equation}
on $\calK(\lambda).$ Introducing
$\calH(\lambda)=\overline{(\calH;
(\cdot,\cdot)_{\lambda})},$ the completion of $\calH$ with
respect to the semi-inner product
\begin{equation}
(f,g)_\lambda=(f,A(\lambda)g)_\calH \text{ for a.e. }
\lambda\in (a,b), \lb{4.43f}
\end{equation}
the isometric isomorphism between $\calK(\lambda)$ and
$\calH(\lambda),$ $\lambda\in (a,b)$ then yields that
\eqref{4.43e} is unitarily equivalent to (cf.~\cite{Ka78}), 
\begin{equation}
S(\lambda)=\Phi(\lambda +i0)^{-1}\Phi(\lambda-i0) \text{
for a.e. } \lambda\in (a,b). \lb{4.43g}
\end{equation}
Arguing as in section~5 of Carey~\cite{Ca76}, Asano's
result~\cite{As67} on strong boundary values for
vector-valued singular integrals of Cauchy-type then yields
the following connection between $S(\lambda)$ in
\eqref{4.43g} and $\Xi(\lambda),$ $\lambda\in (a,b)$ in
\eqref{4.35}.

\begin{lemma} \lb{l4.7} Assume Hypothesis~\ref{h3.1} with
$J=I_\calH$ and suppose $H_0$ is spectrally absolutely
continuous on $(a,b)\subseteq\bbR.$ Then $S(\lambda)$ given
by \eqref{4.43g} satisfies
\begin{align}
S(\lambda)&=\exp{\bigg(-\text{P.V.}\int_\bbR d\mu\,\Xi(\mu)
(\mu-\lambda)^{-1} -i\pi\Xi(\lambda)\bigg)}\times \no \\
&\times \exp{\bigg(\text{P.V.}\int_\bbR d\mu\,\Xi(\mu)
(\mu-\lambda)^{-1} -i\pi\Xi(\lambda)\bigg)} \text{ for
a.e. } \lambda\in (a,b), \lb{4.43h}
\end{align}
where $\text{P.V.}\int_\bbR d\mu \, \cdot$ denotes the
principal value. This implies the 
Birman-Krein formula 
\begin{equation}
\text{det}_{\calH (\lambda)} (S(\lambda))=e^{-2\pi
i\xi(\lambda)}
\text{ for a.e. } \lambda\in (a,b), \lb{4.43i}
\end{equation}
with $\xi(\lambda)=\tr_\calH(\Xi(\lambda)).$
\end{lemma}
\begin{proof}
Asano's result \cite{As67}, applied to the Hilbert space of
$\calB_2(\calH)$-operators yields
\begin{align}
&\slim_{\varepsilon\downarrow 0}\int_\bbR
d\mu\,\Xi(\mu)(\mu-(\lambda\pm
i\varepsilon))^{-1} \no \\
&=\text{P.V.}\int_\bbR
d\mu\,\Xi(\mu)(\mu-\lambda)^{-1} \pm
i\pi\Xi(\lambda)\in\calB_2(\calH) 
\text{ for a.e. } \lambda\in\bbR \lb{4.43j} 
\end{align}
(with $\slim$ denoting convergence in
$\calB_2(\calH)$-norm). Combining \eqref{4.43j}
\eqref{4.43g}, and \eqref{4.35} then yields \eqref{4.43h}.
\end{proof}

It should be noted that \eqref{4.43g} is not necessarily
the usually employed scattering operator at fixed energy
$\lambda\in (a,b).$ In concrete applications one infers
that $0\leq A(\lambda)\in\calB_1(\calH)$ typically factors
into a product 
\begin{equation}
A(\lambda)=B(\lambda)^*B(\lambda), \quad
B(\lambda)\in\calB_2(\calH,\calL), \lb{4.43k}
\end{equation}
with $\calL$ another Hilbert space (e.g.,
$\calL=L^2(S^{n-1})$ in connection with potential
scattering in $\bbR^n,$ $n\geq 2$ and $\calL=\bbC^2$ for
$n=1$) and hence usually $S(\lambda)$ in \eqref{4.43g} is
then repaced by the unitary operator
\begin{equation}
\widetilde S(\lambda)=I_\calL-2\pi
iB(\lambda)(I_\calH+\Phi(\lambda+i0))^{-1}B(\lambda)^*
\text{ for a.e. } \lambda\in (a,b) \lb{4.43l}
\end{equation}
in $\calL.$ 

For brevity we only considered positive perturbations
$V=KK^*.$ The general case of perturbations of the type
$V=KJK^*$ will be considered elsewhere \cite{GM99}.

Next, assume Hypothesis~\ref{h3.1} and $J=I_\calH$,
introduce
\begin{equation}\lb{4.44}
V=KK^*,
\end{equation}
and define the family of self-adjoint operators
\begin{equation}\lb{4.45}
H(s)=H_0+sV, \quad s\in \bbR.
\end{equation}
In accordance with \eqref{3.5}--\eqref{3.7}
introduce the following bounded operators
\begin{align}
\Phi_+(z, s)
&=I_{\calH}+ sK^*(H_0-z)^{-1}K, \quad s>0,
\lb{4.46} \\
\widetilde \Phi_-(z,s)
&=I_{\calH}+s K^*(H_+-z)^{-1}K, \quad s<0,\lb{4.46b}
\end{align}
for $z\in \bbC\backslash \bbR$.

By Lemma~\ref{l3.4} we have the representations
\begin{align}
&\log(\Phi_+(z,s))=\int_{\bbR}d\lambda \,
\Xi_+(\lambda,s)(\lambda-z)^{-1}, \quad s>0
\lb{4.47} \\
&\log(\widetilde \Phi_-(z,s))=-\int_{\bbR}d\lambda 
\, \Xi_-
(\lambda,s)(\lambda-z)^{-1}, \quad s<0, \lb{4.47b}
\end{align}
where $\Xi_+(\lambda,s)$ (respectively,  $\Xi_-(\lambda,s)$) 
is the spectral shift operator associated with the pair
$(H_0, H(s))$ for $s>0$ (respectively, for $s<0$).

Moreover, by \eqref{3.26aa} and \eqref{3.26ab} we have
\begin{align}
\Xi_+(\lambda,s)&=\lim_{\varepsilon\downarrow 0}
\pi^{-1}\Im(\log(\Phi_+(\lambda+i\varepsilon, s))), 
\quad s>0, \lb{4.48} \\
\Xi_-(\lambda,s)&=-\lim_{\varepsilon\downarrow 0} 
\pi^{-1}\Im(\log(\wti \Phi_-(\lambda+i\varepsilon,s))), 
\quad s<0. \lb{4.49}
\end{align}

\begin{theorem}\lb{t4.8}
Assume Hypothesis~\ref{h3.1} and $J=I_\calH.$ Set 
$V=KK^*\ge 0,$ suppose that $\calP=\ran(V)$ 
is a finite-dimensional subspace of $\calH,$ and 
denote by $P$ the orthogonal projection onto $\calP$. 
In addition, let
\begin{equation}\lb{4.50}
T(z)=V^{1/2}(H_0-z)^{-1}V^{1/2}, \quad z\in \bbC_+.
\end{equation}
Then the boundary values
\begin{equation} \lb{4.51}
T(\lambda)=\nlim_{\varepsilon\downarrow 0} 
T(\lambda+i\varepsilon)
\end{equation}
exist for a.e.~$\lambda\in \bbR.$ For such 
$\lambda\in\bbR,$ 
$T(\lambda)$ is reduced by the subspace $\calP$ 
and the part $T(\lambda)|_\calP$ of the operator 
$T(\lambda)$ restricted to the subspace $\calP$ is 
invertible for a.~e.~$\lambda\in\bbR.$ The 
corresponding set of $\lambda\in\bbR$ such that 
$T(\lambda)|_\calP$ is invertible in $\calP=P\calH$ 
is denoted by $\Lambda.$ Finally, for all 
$\lambda\in \Lambda$ one obtains the
following  asymptotic expansion
\begin{equation}\lb{4.52}
\Xi_+(\lambda, s)+\Xi_-(\lambda,- s)
\underset{s\uparrow\infty}{=}
P-2(\pi s)^{-1} P\, \Im (T(\lambda)|_\calP)^{-1}P+ 
O(s^{-2})P.
\end{equation}
\end{theorem}
\begin{proof} The a.e.~existence of the norm limit in 
\eqref{4.51} and the invertibility of $T(\lambda)|_\calP$ 
in $\calP=P\calH$ is a consequence of Lemma~\ref{l2.3b}. 
By  definition \eqref{2.17} of logarithms of dissipative  
operators one infers
\begin{align}
\log(\Phi_+(\lambda+i\varepsilon, s))=
-i \int_0^\infty
dt\big( (sT(\lambda+i\varepsilon)+(1+it) 
I_\calH)^{-1}-(1&+it)^{-1}
I_\calH\big),
\no \\ 
& s>0, \, \varepsilon >0.\lb{4.53}
\end{align}
By \eqref{4.50}, $\log(\Phi_+(\lambda+i\varepsilon, s))$ 
is reduced by the subspace $\calP=P\calH$ and 
\begin{equation}\lb{4.54}
\log(\Phi_+(\lambda+i\varepsilon, 
s))|_{\calH\ominus \calP}=0.
\end{equation}
The operator $\log(\Phi_+(\lambda+i\varepsilon, 
s))|_{ \calP}$ restricted to the invariant subspace
$\calP$ then can be represented as follows
\begin{align}
&\log(\Phi_+(\lambda+i\varepsilon, s))|_{ \calP}
\no \\
&=-i \int_0^\infty
dt\big( (sT(\lambda+i\varepsilon)|_{\calP}
+(1+it) I_{\calP})^{-1}-(1+it)^{-1}
I_{\calP}\big), \lb{4.55} \\ 
&\hspace{7.35cm} s>0, \, \varepsilon >0. \no
\end{align}
For $\lambda\in\Lambda,$ the operator 
$(I+sT(\lambda))|_{\calP}$ is  invertible for $s>0$ 
sufficiently large and therefore, for such $s>0$ one can go 
to the limit $\varepsilon \to 0$ in \eqref{4.55} to 
arrive at
\begin{align}
&\log(\Phi_+(\lambda+i0, s))|_{ \calP} \no \\
&=-i \int_0^\infty
dt\big( (sT(\lambda)|_{\calP}+(1+it) I_{\calP})^{-1}
-(1+it)^{-1} I_{\calP}\big), \lb{4.56} \\ 
& \hspace{4cm} s>0 \text{ sufficiently large, } 
\lambda\in\Lambda. \no
\end{align}

Since for $s<0$ the operator $-\wti 
\Phi_-(\lambda+i\varepsilon,s)$ is 
dissipative, one concludes as in \eqref{2.55a} that
\begin{align}
&\log(\widetilde \Phi_-(\lambda+i\varepsilon,s))=
\big (\log(\widetilde \Phi_-^*(\lambda
+i\varepsilon,s))\big )^* \no \\
&=i\int_0^\infty
dt\big( (sT(\lambda+i\varepsilon)+(1-it) 
I_\calH)^{-1}-(1-it)^{-1}
I_\calH\big) \no \\
&=-i \int_0^\infty
dt\big( (|s|T(\lambda+i\varepsilon)+(it-1) 
I_\calH)^{-1}+(1-it)^{-1} I_\calH\big), \lb{4.57} \\
&\hspace{7.25cm} s<0, \, \varepsilon >0. \no
\end{align}

Similarly one concludes that  
$\log(\widetilde \Phi_-(\lambda+i0,s)),$ 
$\lambda\in\Lambda$ is reduced 
by the subspace $\calP$ and
\begin{align}
&\log(\widetilde \Phi_-(\lambda+i0,s))|_{ \calP}
\no \\
&=-i \int_0^\infty
dt\big( (sT(\lambda)|_{\calP}+(it-1) I_{\calP})^{-1}
-(1+it)^{-1} I_{\calP}\big), \lb{4.58} \\ 
& \hspace{4cm} s>0 \text{ sufficiently large, }
\lambda\in\Lambda. \no
\end{align}

By \eqref{4.48} and \eqref{4.49} one obtains for $s>0$ 
\begin{align} 
&\Xi_+(\lambda, s)+\Xi_-(\lambda,- s)=
\pi^{-1}\Im \big (
\log(\Phi_+(\lambda+i0, s))-
\log(\widetilde \Phi_-(\lambda+i0,s)) \big ), 
\lb{4.59} \\ 
& \hspace{10.2cm} \lambda\in\Lambda. \no
\end{align}
Combining \eqref{4.56} and \eqref{4.58} and taking into 
account the fact that
\begin{equation}\lb{4.60}
\log(\Phi_+(\lambda+i0, s))|_{\calH\ominus \calP}=
\log(\widetilde \Phi_-(\lambda+i0,s))|_{\calH\ominus 
\calP}=0, \quad \lambda\in\Lambda,
\end{equation}
one concludes that the subspace $\calP$ reduces 
$\Xi_+(\lambda, s)+\Xi_-(\lambda,- s)$ and that
\begin{equation}\lb{4.61}
\big (\Xi_+(\lambda, s)+\Xi_-(\lambda,- 
s)\big )|_{\calH\ominus \calP}=0, \quad 
\lambda\in\Lambda.
\end{equation}
Moreover, for $\lambda\in\Lambda,$
\begin{align}
&\big (\Xi_+(\lambda, s)+\Xi_-(\lambda,- 
s)\big )|_{ \calP} \no \\
&=\pi^{-1}\Im \bigg (-i \int_0^\infty
dt\big( (sT(\lambda)|_{\calP}+(1+it) I_{\calP})^{-1}
-(sT(\lambda)|_{\calP}+(it-1) I_{\calP})^{-1}
\bigg ) \no \\
&+\pi^{-1} \Im \bigg (i \int_0^\infty dt 
\big (((1+it)^{-1}+(1-it)^{-1})I_{\calP} \bigg ) \no \\ 
&=\pi^{-1} \Im \bigg (2i \int_0^\infty
dt\big( (sT(\lambda)|_{\calP}+(it+1) I_{\calP})^{-1} 
(sT(\lambda)|_{\calP} +(it-1) I_{\calP})^{-1}\big) 
\bigg )  +I|_{\calP}. \no \\
\lb{4.62}
\end{align}
Changing variables $t\to s^{-1}t$ using the fact that
 $T(\lambda)|_{\calP}$ is invertible for 
$\lambda\in\Lambda$ then yields 
\begin{align}
&\int_0^\infty
dt\big( (sT(\lambda)|_{\calP}+(it+1) I_{\calP})^{-1}
(sT(\lambda)|_{\calP}+(it-1) I_{\calP})^{-1}\big)
\no \\
&=s^{-1}\int_0^\infty
dt\big( (T(\lambda)|_{\calP}+(it+s^{-1}) I_{\calP})^{-1}
(T(\lambda)|_{\calP}+(it-s^{-1}) I_{\calP})^{-1}\big)
\no \\
&\underset{s\uparrow\infty}{=}s^{-1}\int_0^\infty
dt\big( (T(\lambda)|_{\calP}+it I_{\calP})^{-2}+O(s^{-2})P
\no \\
&\underset{s\uparrow\infty}{=}s^{-1}i 
((T(\lambda)|_{\calP})^{-1}+O(s^{-2})P, \lb{4.63}
\end{align}
where in obvious notation $O(s^{-2})$ denotes a bounded 
operator in $\calP$ whose norm is of order $O(s^{-2})$ 
as $s\uparrow\infty.$ Combining
\eqref{4.62} and
\eqref{4.63} we get the  asymptotic representation
\begin{equation}\lb{4.64}
\big (\Xi_+(\lambda, s)+\Xi_-(\lambda,-s) 
\big )|_{ \calP}\underset{s\uparrow\infty}{=}
I|_\calP -2(\pi s)^{-1} \Im 
\big ((T(\lambda)|_{\calP})^{-1}\big ) +O(s^{-2})P, 
\quad \lambda\in\Lambda.
\end{equation}
Together with \eqref{4.61} this proves \eqref{4.52}.
\end{proof}

Taking the trace of \eqref{4.52} and going 
to the limit $s\uparrow\infty,$ Theorem~\ref{t4.8}
implies the following result first proved by Simon (see
\cite{Si98}) for finite-rank nonnegative perturbations $V,$
\begin{equation}\lb{4.65}
\lim_{s\uparrow\infty} (\xi(\lambda, H_0,H_0+sV)
-\xi(\lambda, H_0,H_0-sV))=\text{rank}(V), \quad 
\lambda\in\Lambda,
\end{equation}
where
\begin{equation}
\xi(\lambda, H_0,H_0+sV)=\tr (\Xi_+(\lambda, s)), 
\quad s>0 \lb{4.66}
\end{equation}
and
\begin{equation}
\xi(\lambda, H_0,H_0-sV)=-\tr( \Xi_-(\lambda, s)), 
\quad s>0 \lb{4.67}
\end{equation}
are spectral shift functions associated with the pairs
$(H_0, H+sV),$  and  $(H_0, H-sV),$ $ s>0,$ respectively.

Finally we turn to an application concerning an approach
to abstract trace formulas based on perturbation 
theory for 
pairs of self-adjoint extensions  of a common closed, 
symmetric,
densely defined linear  operator $\dot H$ in some complex
separable Hilbert space
$\calH.$ We first  treat the simplest case of deficiency 
indices
$(1,1)$ and hint at extensions to the case of deficiency
indices $(n,n),$ $n\in\bbN$ at the end. These results are 
 applicable to  one-dimensional (matrix-valued) 
Schr\"odinger
operators.

We start by setting up the basic formalism. Assuming 
\begin{equation}
\text{def}\,(\dot H) =(1,1), \lb{4.68}
\end{equation} 
we use von Neumann's parametrization of all self-adjoint 
extensions $H_\alpha,$ $\alpha \in [0,\pi),$ in the 
usual form
\begin{align}
&H_\alpha (f+c(u_+ +e^{2i\alpha}u_-))=\dot H f+c(iu_+  
-ie^{2i\alpha}u_-), \quad \alpha\in [0,\pi), \no \\
&\dom(H_\alpha)=\{f+c(iu_+ - ie^{2i\alpha}u_-)\in 
\dom({\dot H}^*) \,|\, f\in\dom(\dot H), \, c\in\bbC \}, 
 \lb{4.69}
\end{align}
where 
\begin{equation}
u_\pm \in \dom({\dot H}^*), \quad {\dot H}^*u_\pm 
=\pm iu_\pm, \quad \|u_\pm\|_\calH=1.
\end{equation}
Introducing the Donoghue $m$-function (cf.~\cite{Do65}, 
\cite{GKMT98}, \cite{GMT98}, \cite{GT97}) associated with 
the pair $(\dot H,H_\alpha)$ by
\begin{equation}
m_{\alpha}(z)=z+(1+z^2)(u_+,(H_{\alpha}-z)^{-1}u_+)_{\calH}, 
\quad z\in\bbC\backslash\bbR, \, \alpha \in [0,\pi), 
\lb{4.70}
\end{equation}
one verifies 
\begin{equation}
m_{\beta}(z)=\frac{-\sin(\beta-\alpha) +
\cos(\beta-\alpha) m_{\alpha}(z)}{\cos(\beta-\alpha) +
\sin(\beta-\alpha) m_{\alpha}(z)}, \quad
\alpha,\beta\in [0,\pi), \lb{4.71}
\end{equation}
and obtains Krein's formula for the resolvent difference of 
two self-adjoint extensions $H_\alpha,$ $H_\beta$ of 
$\dot H$ (cf.,~\cite{AG93}, Sect.~84, \cite{GMT98}),
\begin{align}
&(H_\alpha-z)^{-1}-(H_\beta-z)^{-1}=(m_\alpha(z)+ 
\cot(\beta-\alpha))^{-1}(u_+(\overline z),\cdot)_\calH u_+(z), 
\no \\
&\hspace*{7.3cm} z\in\bbC\backslash\bbR, \, \alpha,
\beta\in [0,\pi), \lb{4.72}
\end{align} 
where
\begin{equation}
u_+(z)=(H_\alpha -i)(H_\alpha -z)^{-1}, \quad
z\in\bbC\backslash\bbR. \lb{4.73}
\end{equation}
For later reference we note the useful facts,
\begin{align}
&(u_+(\overline z_1),u_+(z_2))_\calH=\frac{m_\alpha(z_1)-
m_\alpha(z_2)}{z_1-z_2}, \quad z_1,z_2 \in\bbC\backslash\bbR, 
\lb{4.74} \\
&(u_+(\overline z),u_+(z))_\calH= \frac{d}{dz}m_\alpha (z), 
\quad z \in\bbC\backslash\bbR. \lb{4.75}
\end{align}
Next, we consider a bounded self-adjoint operator $V$ in 
$\calH,$
\begin{equation}
V=V^*\in\calB(\calH) \lb{4.76}
\end{equation}
and introduce an ``unperturbed'' operator
${\dot H}^{(0)}=\dot H - V$  in $\calH$ such that 
\begin{equation}
\dot H = {\dot H}^{(0)} + V, \quad \dom(\dot H) =
\dom({\dot H}^{(0)}). \lb{4.77}
\end{equation}
Consequently, 
\begin{equation}
{\dot H}^* = {\dot H}^{(0)*} + V, \quad \dom({\dot H}^*) =
\dom({\dot H}^{(0)*}). \lb{4.77a}
\end{equation}
In addition, we pick $\alpha, \alpha^{(0)}\in [0,\pi)$ 
such that 
\begin{equation}
\dom(H_\alpha)=\dom(H_{\alpha^{(0)}}^{(0)}). \lb{4.77b}
\end{equation}
Formulas \eqref{4.68}--\eqref{4.75} then apply to the
self-adjoint  extensions of ${\dot H}^{(0)}$ and in obvious
notation we denote  corresponding quantities associated with
${\dot H}^{(0)}$ by ${\dot H}^{(0)*},$ 
$H^{(0)}_{\alpha^{(0)}},$ $u^{(0)}_\pm,$ 
$u^{(0)}_+ (z),$ $m_{\alpha^{(0)}}^{(0)}(z),$ 
$\alpha^{(0)}\in [0,\pi),$ etc.  A fundamental link between
${\dot H}^{(0)*}$ and ${\dot H}^*$ is  provided by the 
following
result.
\begin{lemma} \lb{l4.9}
Assume \eqref{4.76} and \eqref{4.77} and let $z\in
\bbC\backslash\bbR.$ Then
$(I_\calH-(H_\alpha  -z)^{-1}V),$ is
invertible,
\begin{equation}
(I_\calH-(H_\alpha -z)^{-1}V)^{-1}=
(I_\calH+(H_{\alpha^{(0)}}^{(0)} -z)^{-1}V) \lb{4.78}
\end{equation}
and 
\begin{equation}
\ker({\dot H}^* -zI_\calH)=(I_\calH-(H_\alpha -z)^{-1}V)
\ker({\dot H}^{(0)*} -zI_\calH).
\lb{4.79}
\end{equation}
In particular,
\begin{equation}
u_+(z)=c(I_\calH-(H_\alpha-z)^{-1}V)u_+^{(0)}(z), \lb{4.80}
\end{equation}
where $c>0$ is determined by the requirement
$\|u_+(i)\|_\calH=1.$
\end{lemma}
\begin{proof}
Equation \eqref{4.78} is clear from the identities
\begin{align}
I_\calH&=(I_\calH-(H_\alpha -z)^{-1}V)
(I_\calH+(H_\alpha^{(0)}
-z)^{-1}V) \no \\
&=(I_\calH+(H_\alpha^{(0)} -z)^{-1}V)(I_\calH-(H_\alpha
-z)^{-1}V) \lb{4.81}
\end{align}
and \eqref{4.71} follows since  
$(H_\alpha -z)^{-1}$ maps $\calH$ into $\dom({\dot H}^*)= 
\dom({\dot H}^{(0)*})$ and hence 
\begin{align}
&({\dot H}^* -z)(I_\calH-(H_\alpha-z)^{-1}V)g =
({\dot H}^* -z)g-Vg \no \\
&=({\dot H}^{(0)*}+V-z)g-Vg=({\dot H}^{(0)*}-z)g=0, 
\quad g\in\ker({\dot H}^{(0)*}-zI_\calH). \lb{4.82}
\end{align}
Equation~\eqref{4.80} is then clear from \eqref{4.78}, 
\eqref{4.79}.
\end{proof}

In the following we assume in addition that 
${\dot H}^{(0)}$ 
is bounded from below, that is,
\begin{equation}
{\dot H}^{(0)}\geq CI_\calH \text{ for some } C\in\bbR
\lb{4.83}
\end{equation}
and choose $\beta,\beta^{(0)}\in [0,\pi)$ such that 
\begin{equation}
\dom(H_{\beta})=\dom(H_{\beta^{(0)}}^{(0)})  \lb{4.84}
\end{equation}
and denote the Friedrichs extensions of $\dot H$ and 
${\dot H}^{(0)}$ by $H_{\alpha_F}$ and 
$H_{\alpha_F^{(0)}}^{(0)},$ respectively. In particular, 
since $V\in\calB(\calH)$ this implies 
\begin{equation}
\dom(H_{\alpha_F})=\dom(H_{\alpha^{(0)}_F}^{(0)}).
\lb{4.85}
\end{equation}
Throughout the remainder of this section, the 
subscript $F$
indicates the Friedrichs extension of ${\dot H}^{(0)}$ 
and $\dot H$ and we choose $\alpha=\alpha_F$ ($\alpha^{(0)}
=\alpha^{(0)}_F$) in \eqref{4.74}, \eqref{4.81}, etc.
We recall  (cf.,~e.g., \cite{Do65}, \cite{GT97}) that 
$\alpha_F$
for  the Friedrichs extension $H_{\alpha_F}$ of 
$\dot H$ (and
similarly
$\alpha^{(0)}_F$ for the Friedrichs extension 
$H^{(0)}_{\alpha^{(0)}_F}$ of ${\dot
H}^{(0)}$)  is uniquely characterized by 
\begin{equation}
\lim_{z\downarrow -\infty} m_{\alpha_F}(z)=-\infty. 
\lb{4.86}
\end{equation}
Next, taking into account \eqref{4.72}, \eqref{4.84}, and
\eqref{4.85}, we recall the exponential Herglotz 
representations (cf.~also \cite{KJ81}),
\begin{equation}
\ln(m_{\alpha_F}(z)+\cot(\beta-\alpha_F))=
c_{\alpha_F,\beta}+ \int_\bbR d\lambda\,
((\lambda-z)^{-1}-\lambda(\lambda^2
+1)^{-1})\eta_{\alpha_F,\beta}(\lambda), \lb{4.93} 
\end{equation}
\begin{align}
&\ln(m_{\alpha_F^{(0)}}^{(0)}(z)+\cot(\beta^{(0)}-
\alpha_F^{(0)})) \lb{4.94} \\
&= c_{\alpha_F^{(0)},\beta^{(0)}}
 +\int_\bbR d\lambda\, ((\lambda-z)^{-1}-\lambda(\lambda^2
+1)^{-1})\eta_{\alpha_F^{(0)},\beta^{(0)}}^{(0)}(\lambda), 
\no \\  
& \hspace*{.8cm} 0\leq \eta_{\alpha_F,\beta}(\lambda),
\eta_{\alpha_F^{(0)},\beta^{(0)}}^{(0)}(\lambda)\leq 1  
\text{ for a.e. } \lambda\in\bbR. \lb{4.95}
\end{align}

Combining the paragraph following \eqref{3.20} with 
\eqref{3.21}, \eqref{4.72}, \eqref{4.75} and
\eqref{4.93}, one verifies that 
$\eta_{\alpha_F,\beta}(\lambda)$ represents the Krein 
spectral
shift function for the  pair $(H_{\alpha_F},H_\beta)$ (and
analogously in the  unperturbed case).

The following result links $m_{\alpha_F} (z)$ and 
$m_{\alpha_F^{(0)}}^{(0)}(z).$
\begin{lemma} \lb{l4.10}
Suppose $z<0,$ $|z|$ sufficiently large. Then
\begin{align}
&\frac{d}{dz}m_{\alpha_F}(z)\underset{z\downarrow 
-\infty}{=}
c^2\frac{d}{dz}m^{(0)}_{\alpha_F^{(0)}}(z) -c^2\frac{d}{dz}
(u^{(0)}_+(z),Vu^{(0)}_+(z))_\calH + O(|z|^{-2}), 
\lb{4.87} \\
&m_{\alpha_F}(z)\underset{z\downarrow -\infty}{=}
c^2m^{(0)}_{\alpha_F^{(0)}}(z) +C_F -c^2
(u^{(0)}_+(z),Vu^{(0)}_+(z))_\calH + O(|z|^{-1}), \lb{4.88}
\end{align}
where 
\begin{equation}
C_F=c^2\cot(\beta^{(0)}-\alpha_F^{(0)})-
\cot(\beta-\alpha_F). \lb{4.89} 
\end{equation} 
\end{lemma}
\begin{proof}
Using \eqref{4.75}, \eqref{4.80}, 
\begin{equation}
(H_{\alpha_F}-z)^{-1}u_+(z)=(H_{\alpha_F}-i)
(H_{\alpha_F}-z)^{-2}=\frac{d}{dz}u_+(z), \lb{4.90}
\end{equation}
the resolvent equation
\begin{equation}
(H_{\alpha_F}-z)^{-1}=(H_{\alpha_F^{(0)}}^{(0)}-z)^{-1} 
-(H_{\alpha_F^{(0)}}^{(0)}-z)^{-1}V(H_{\alpha_F}-z)^{-1}, 
\quad z\in\bbC\backslash\bbR, \lb{4.91}
\end{equation}
and 
\begin{equation}
\|u_+^{(0)}(z)\|_\calH\underset{z\downarrow -\infty}{=}O(1), 
\lb{4.92}
\end{equation}
one verifies \eqref{4.87}. Equation \eqref{4.88} then 
follows 
upon integrating \eqref{4.87}. The integration constant
$C_F$  can be determined by a somewhat lengthy perturbation
argument as follows. For brevity we will temporarily use the
following short-hand notations,
\begin{align}
&M(z)=m_{\alpha_F}(z), \quad m(z)=
m^{(0)}_{\alpha^{(0)}_F}(z),
\quad  v(z)=(u^{(0)}_+(z),Vu^{(0)}_+(z))_\calH, \no \\
&\gamma=\cot(\beta-\alpha_F), \quad
\delta=\cot(\beta^{(0)}-\alpha_F^{(0)}).
\lb{4.92a}
\end{align}
First we claim
\begin{equation}
\lim_{z\downarrow -\infty}z^{-2}m(z)^2/m'(z)=0. \lb{4.92b}
\end{equation}
Since none of the spectral measures $d\mu_\beta(\lambda),$
$\beta\in [0,\pi)$ associated  with the Herglotz
representation  of $m_\beta(z)$ in \eqref{4.71} is a finite
measure on $\bbR,$ one obtains from
$m'(z)/(m(z)+\delta)^2=((-(m(z)+\delta)^{-1}))'$ 
\begin{equation}
z^2m'(z)(m(z)+\delta)^{-2}=\int_\bbR
d\mu_{\beta_0}(\lambda)\,z^2(\lambda -z)^{-2} \uparrow
+\infty \text{ as } z\downarrow  -\infty \lb{4.92d}
\end{equation}
for some $\beta_0\in [0,\pi).$ 

Next we will show that
\begin{align}
&(d/dz)\ln((M(z)+\gamma)/(m(z)+\delta)) \lb{4.92e} \\ 
&=\int_\bbR d\lambda\,
(\lambda-z)^{-2}(\eta_{\alpha_F,\beta}(\lambda)-
\eta_{\alpha_F^{(0)},\beta^{(0)}}^{(0)}(\lambda)) 
\lb{4.92f} \\
&=\tr((H_{\alpha_F}-z)^{-1}-(H_\beta-z)^{-1} -
(H^{(0)}_{\alpha_F^{(0)}}-z)^{-1}+(H^{(0)}_{\beta^{(0)}}
-z)^{-1}) \lb{4.92g} \\
&\underset{z\downarrow -\infty}{=} O(z^{-2}) \lb{4.92h}.
\end{align}
While \eqref{4.92e}--\eqref{4.92g} are clear from 
\eqref{4.72}, 
\eqref{4.93}, and \eqref{4.94}, we need to prove the 
asymptotic 
relation \eqref{4.92h}. The difference of the first 
and the third resolvent as well as the difference of the 
second and fourth resolvent under the trace in 
\eqref{4.92g} 
is clearly of $O(z^{-2})$ in norm using the resolvent 
equation 
and the fact that $V$ is a bounded operator. On the 
other hand, 
the operator under the trace in \eqref{4.92g} is the 
difference
of  two rank-one operators by \eqref{4.72} and hence 
at most of
rank two. Hence the trace norm of the operator under 
the trace
in \eqref{4.92g} is also of order $O(z^{-2})$ as 
$z\downarrow 
-\infty.$ 

Integrating \eqref{4.92e} taking into account 
\eqref{4.92h} then proves that 
\begin{equation}
m(z)/M(z)=O(1) \text{ and } M(z)/m(z)=O(1) \text{ as } 
z\downarrow -\infty. \lb{4.92i}
\end{equation}

Next we abbreviate
\begin{equation}
D(z)=(d/dz)\ln((M(z)+\gamma)/(m(z)+\delta))
\end{equation}
and compute (cf.~\eqref{4.87})
\begin{align}
D(z)&=\frac{M'(z)(m(z)+\delta)-m'(z)(M(z)+\gamma)}
{(M(z)+\gamma)(m(z)+\delta)} \no \\
&=\frac{(c^2m'(z)+r(z))(m(z)+\delta)-m'(z)(M(z)+\gamma)}
{(M(z)+\gamma)(m(z)+\delta)} \no \\
&=\frac{m'(z)}{(M(z)+\gamma)(m(z)+\delta)}
(c^2(m(z)+\delta)-
(M(z)+\gamma)) + \frac{r(z)}{M(z)+\gamma}, \lb{4.92j}
\end{align}
where
\begin{equation}
r(z)=-v'(z)+O(z^{-2})=o(|m'(z)/m(z)|) 
\text{ as } z\downarrow -\infty  \lb{4.92k}
\end{equation}
by taking into account
\begin{equation}
m(z)=o(|z|) \text{ as } z\downarrow -\infty  \lb{4.92ka}
\end{equation} 
and
\begin{equation}
v'(z)=O(|m'(z)/z|) \text{ as } z\downarrow -\infty, 
\lb{4.92l}
\end{equation}
which in turn follows from \eqref{4.75}, \eqref{4.90}, 
the
fact  that $V$ is a bounded operator, and
$\|(H^{(0)}_{\alpha_F^{(0)}}-z)^{-1}\|=O(|z|^{-1})$ as 
$z\downarrow -\infty.$ Thus,
\begin{align}
&(m(z)^2/m'(z))D(z) \no \\
&=\frac{m(z)^2}{(M(z)+\gamma)(m(z)+\delta)}
(c^2(m(z)+\delta)
-(M(z)+\gamma)) \no \\
&\hspace*{.4cm} + \frac{m(z)^2r(z)}{m'(z)(M(z)+\gamma)} 
\no \\
&=\frac{m(z)^2}{(M(z)+\gamma)(m(z)+\delta)}
(c^2(m(z)+\delta)
-(M(z)+\gamma)) + o(1) \no \\ 
&=-(m(z)/M(z))((M(z)+\gamma)-c^2(m(z)+\delta))) +o(1) 
\no \\
&=o(1) \text{ as } z\downarrow -\infty, \lb{4.92m}
\end{align}
by \eqref{4.92b} and \eqref{4.92h}. This implies
\begin{equation}
\lim_{z\downarrow -\infty}((M(z)+\gamma)-
c^2(m(z)+\delta))=0
\end{equation}
by \eqref{4.92i} and hence proves \eqref{4.88} and 
\eqref{4.89}.
\end{proof} 

Our principal asymptotic result then reads as follows.
\begin{theorem} \lb{t4.11}
\begin{align}
&\int_\bbR d\lambda \, (\lambda-z)^{-2}(\eta_{\alpha_F}
(\lambda)-\eta_{\alpha_F^{(0)}}^{(0)}(\lambda)) \no \\
&\underset{z\downarrow -\infty}{=} 
-\frac{d}{dz}\bigg(\frac{(u^{(0)}_+(z),Vu^{(0)}_+(z))_\calH}
{m^{(0)}_{\alpha^{(0)}_F}(z)
+\cot(\beta^{(0)}-\alpha_F^{(0)})}\bigg) + o(|z|^{-2}) 
\lb{4.96}
\end{align}
and
\begin{equation} 
\underset{z\downarrow -\infty}{\overline{\lim}} \int_\bbR
d\lambda
\, z^2(\lambda-z)^{-2}(\eta_{\alpha_F}
(\lambda)-\eta_{\alpha_F^{(0)}}^{(0)}(\lambda)) 
\text{ exists.} \lb{4.96a}
\end{equation}
\end{theorem}
\begin{proof}
Differentiating \eqref{4.93} and \eqref{4.94} with 
respect to
$z,$ taking into account \eqref{4.87}--\eqref{4.89} yields
\begin{align}
&\int_\bbR d\lambda \, (\lambda-z)^{-2}(\eta_{\alpha_F}
(\lambda)-\eta_{\alpha_F^{(0)}}^{(0)}(\lambda)) \no \\
&\underset{z\downarrow -\infty}{=} 
\frac{(d/dz)m^{(0)}_{\alpha^{(0)}_F}(z)-
(d/dz)(u^{(0)}_+(z),Vu^{(0)}_+(z))_\calH +O(|z|^{-2})}
{m^{(0)}_{\alpha^{(0)}_F}(z)+\cot(\beta^{(0)}-
\alpha_F^{(0)})-
(u^{(0)}_+(z),Vu^{(0)}_+(z))_\calH +O(|z|^{-1})} \no \\
&\hspace*{1cm}-\frac{(d/dz)m^{(0)}_{\alpha^{(0)}_F}(z)}
{m^{(0)}_{\alpha^{(0)}_F}(z)+\cot(\beta^{(0)}-
\alpha_F^{(0)})}. 
\lb{4.97}
\end{align}
In order to verify \eqref{4.96} we need to estimate various
terms. For brevity we will again temporarily use the 
short-hand notations introduced in \eqref{4.92a}. Thus, 
\eqref{4.97}
becomes
\begin{align}
&\frac{m'(z)-v'(z)+O(z^{-2})}{m(z)+\delta-v(z)+O(|z|^{-1})} 
-\frac{m'(z)}{m(z)+\delta} \no \\
&=-(d/dz)(v(z)/(m(z)+\delta)) +O(|m(z)^{-1}z^{-2}|) \no \\
&\hspace*{.42cm} +O(|z^{-1}m'(z)m(z)^{-2}|) +
O(|m'(z)v(z)^2m(z)^{-3}|) \no \\ 
&=-(d/dz)(v(z)/(m(z)+\delta))
+o(z^{-2}) \text{ as }  z\downarrow -\infty \lb{4.97a}
\end{align}
and we need to verify the last line in \eqref{4.97a} 
and the 
claim \eqref{4.96a}. By \eqref{4.86}, one concludes
\begin{equation}
O(|m(z)^{-1}z^{-2}|)=o(z^{-2}) \text{ as }  
z\downarrow -\infty.
\lb{4.97b}
\end{equation}
Next, using 
\begin{equation}
v(z)=O(|m'(z)|) \text{ as }  z\downarrow -\infty \lb{4.97c}
\end{equation}
(cf. \eqref{4.75}), one obtains
\begin{equation}
O(|m'(z)v(z)^2m(z)^{-3}|)=O(|m'(z)^3m(z)^{-3}|)=O(|z|^{-3}) 
\text{ as }  z\downarrow -\infty \lb{4.97d}
\end{equation}
since 
\begin{equation}
m'(z)m(z)^{-1}=O(|z|^{-1}) \text{ as }  z\downarrow -\infty.
\lb{4.97e}
\end{equation}
Relation \eqref{4.97e} is shown as follows. Since $-m'(z)
(m(z)+\delta)^{-1}=
(\ln(-(m(z)+\delta)^{-1})',$ and $-(m(z)+\delta)^{-1}$ 
(being
distinct from the  Friedrichs $m$-function) belongs to some
measure $d\mu_{\beta_0}(\lambda)$ in \eqref{4.71} with
$\smallint_\bbR d\mu_{\beta_0}(\lambda)(1+|\lambda|)^{-1}
<\infty,$ one concludes from \cite{AD56} that also
$\smallint_\bbR
d\lambda\,\xi_{\beta_0}(\lambda)(1+|\lambda|)^{-1} 
<\infty$ in
the corresponding exponential Herglotz representation of
$-(m(z)+\delta)^{-1}$ (cf.~\eqref{4.94}) and hence
\begin{equation}
-zm'(z)(m(z)+\delta)^{-1}=\int_\bbR 
d\lambda\,(-z)(\lambda 
-z)^{-2}\xi_{\beta_0}(\lambda)=O(1) \text{ as } 
z\downarrow 
-\infty. \lb{4.97f}
\end{equation}
Next, we note
\begin{equation}
O(|z^{-1}m'(z)m(z)^{-2}|)=O(|zm'(z)m(z)^{-1}|)
O(|z^{-2}m(z)^{-1}|)=o(z^{-2}) \text{ as } z\downarrow 
-\infty \lb{4.97g}
\end{equation}
by \eqref{4.86} and \eqref{4.97e} which proves \eqref{4.96}. 

To prove \eqref{4.96a} one estimates
\begin{equation}
v'(z)m(z)^{-1}=O(|z^{-1}m'(z)m(z)^{-1}|)=O(z^{-2}) 
\text{ as } z\downarrow -\infty \lb{4.97h}
\end{equation}
using \eqref{4.92l} and \eqref{4.97e}. Similarly,
\begin{equation}
v(z)m'(z)m(z)^{-2}=O(|m'(z)^2m(z)^{-2}|)=O(z^{-2}) 
\text{ as } z\downarrow -\infty \lb{4.97i}
\end{equation}
by \eqref{4.97e}. This completes the proof.
\end{proof}

Next we will show that the abstract asymptotic result
\eqref{4.97} contains concrete trace formulas for
one-dimensional Schr\"odinger operators first derived in 
\cite{GS96} (see also \cite{Ge95}, \cite{GH95},
\cite{GHS95}, \cite{GHSZ95}), and hence can be viewed as an 
abstract approach to trace formulas.

To make the connection with Schr\"odinger operators on the 
real line we choose $\calH=L^2(\bbR;dx),$ pick a $y\in\bbR,$ 
and identify
$V$ with the real-valued  potential $V(x)$ assuming
\begin{equation}
V\in L^\infty(\bbR;dx)\cap C((y-\varepsilon,y+\varepsilon)) 
\text{ for some } \varepsilon >0. \lb{4.98}
\end{equation}
Similarly, we identify (in obvious notation)
$\dot H,$ ${\dot H}^*,$ $H_{\pi/2},$ $H_{\alpha_F}$ with
\begin{align}
&{\dot H}_y=-d^2/dx^2 + V, \lb{4.99} \\
&\dom({\dot H}_y)=\{g\in L^2(\bbR;dx)\,|\, g,g'\in
AC_{\loc} (\bbR); \, \lim_{\varepsilon\downarrow 0}
g(y\pm\varepsilon)=0; \, g''\in L^2(\bbR;dx) \}, \no \\
&{{\dot H}^*}_y=-d^2/dx^2 + V, \lb{4.100} \\
&\dom({\dot H}^*_y)=\{g\in L^2(\bbR;dx)\,|\, g\in
AC_{\loc} (\bbR), g'\in AC_{\loc} (\bbR\backslash\{y\}); \,
 g''\in L^2(\bbR;dx) \}, \no \\
&H_{y,\pi/2}=-d^2/dx^2 + V, \lb{4.101} \\
&\dom(H_{y,\pi/2})=\{g\in L^2(\bbR;dx)\,|\, g,g'\in
AC_{\loc} (\bbR); \, g''\in L^2(\bbR;dx) \}=H^{2,2}(\bbR), 
\no \\
&{H}_{y,F}=-d^2/dx^2 + V, \lb{4.102} \\
&\dom({H}_{y,F})=\{g\in L^2(\bbR;dx)\,|\, g\in
AC_{\loc} (\bbR), g'\in AC_{\loc} (\bbR\backslash\{y\}); \,
\lim_{\varepsilon\downarrow 0} g(y\pm\varepsilon)=0; \no \\
&\hspace*{9.85cm} g''\in L^2(\bbR;dx) \}, \no 
\end{align}
respectively. (We note that $H_{y,\pi/2}$ is actually
independent of $y\in\bbR.$) The corresponding unperturbed 
operators ${\dot H}^{(0)}_y,$ ${\dot H}^{(0)*}_y,$ 
$H_{y,\pi/2}^{(0)},$ $H_{y,F}^{(0)}$ are then defined as in 
\eqref{4.98}--\eqref{4.102} setting $V(x)=0,$ for all 
$x\in\bbR.$

Denoting by $G(z,x,x')$ the Green's function of
$H_{y,\pi/2},$ that is, 
\begin{equation}
G(z,x,x')=(H_{y,\pi/2}-z)^{-1}(x,x'), \quad 
z\in\bbC\backslash\bbR, \, x,x'\in\bbR, \lb{4.102a}
\end{equation}
and observing that 
\begin{align}
&G^{(0)}(z,x,x')=(H_{y,\pi/2}^{(0)}-z)^{-1}(x,x') 
=i2^{-1/2}z^{-1/2}\exp(iz^{1/2}|x-x'|), \lb{4.102b} \\
&\hspace*{7.75cm} z\in\bbC\backslash\bbR, \, 
x,x'\in\bbR, \no
\end{align}
explicit computations then yield the following 
results for $z<0,$ $|z|$  sufficiently large: 
\begin{align}
&u_+(z,x)=G(z,x,y)/\|G(i,\cdot,y)\|_\calH, \quad 
u_+^{(0)}(z,x)=2^{-5/4}iz^{-1/2}\exp(iz^{1/2}|x-y|), 
\lb{4.103} \\
&\tan(\alpha_F)=-\Re(G(i,y,y))/\Im(G(i,y,y)), \quad 
\alpha_F^{(0)}=3\pi/4, \quad \beta=\beta^{(0)}=\pi/2, 
\lb{4.104} \\ 
&\|u_+(z)\|^2_\calH =\Im(G(z,y,y))/\Im(z), \quad 
\|u_+^{(0)}(z)\|^2_\calH = 2^{-1/2}iz^{-1/2}, \lb{4.105} \\
&(u_+^{(0)}(z),Vu_+^{(0)}(z))_\calH=2^{-1/2}iz^{-1/2}V(y) 
+ o(|z|^{-1/2}), \lb{4.106} \\
&(d/dz)(u_+^{(0)}(z),Vu_+^{(0)}(z))_\calH=
-2^{-3/2}iz^{-3/2}V(y) 
+ o(|z|^{-3/2}), \lb{4.107} \\
&m^{(0)}_{y,F}(z)=i(2z)^{1/2} +1, \quad 
\eta^{(0)}_{\alpha_F^{(0)}}(\lambda)=\begin{cases} 
1/2, & \lambda >0, \\ 
1, & \lambda < 0, \end{cases} \lb{4.108} \\
&(d/dz)(\ln(m^{(0)}_{y,F}(z)- 
(u_+^{(0)}(z),Vu_+^{(0)}(z))_\calH)) -(2z)^{-1} 
\underset{z\downarrow -\infty}{=}2^{-1}V(y)z^{-2} 
+o(|z|^{-2}). \lb{4.109}
\end{align}
Moreover, one computes for $z\in\bbC\backslash\bbR$,
\begin{align}
m_{y,\pi/2}(z)&=z+(1+z^2)(\Im(G(i,y,y)))^{-1}
(G(i,\cdot,y),(H_{y,\pi/2}-z)^{-1}
G(i,\cdot,y))_\calH \no \\
&=(G(z,y,y)-\Re(G(i,y,y)))/\Im(G(i,y,y)), \lb{4.110} \\
m_{y,\alpha_F}(z)&=z+(1+z^2)(\Im(G(i,y,y)))^{-1}
(G(i,\cdot,y),
(H_{y,\alpha_F}-z)^{-1}G(i,\cdot,y))_\calH \no \\
&=(-G(z,y,y)^{-1}|G(i,y,y)|^2+\Re(G(i,y,y)))/\Im(G(i,y,y)).
\lb{4.111} 
\end{align}
Combining \eqref{4.96} and \eqref{4.103}--\eqref{4.108},  
identifying $\eta_{\alpha_F}(\lambda),$ 
$\eta_{\alpha_F^{(0)}}^{(0)}(\lambda)$ with 
$\eta_{\alpha_F}(\lambda,y),$
$\eta_{\alpha_F^{(0)}}^{(0)}(\lambda),$ then  yields 
the trace formula
\begin{equation}
V(y)=\lim_{z\downarrow -\infty}2\int_\bbR d\lambda \,
z^2(\lambda-z)^{-2}(\eta_{\alpha_F}
(\lambda,y)-\eta_{\alpha_F^{(0)}}^{(0)}(\lambda)). 
\lb{4.112}
\end{equation}
Since by \eqref{4.104} and \eqref{4.111},
\begin{equation}
m_{y,\alpha_F}(z) +\tan(\alpha_F)=-G(z,y,y)^{-1}
/\Im(G(i,y,y)), \lb{4.113}
\end{equation}
one can use the exponential Herglotz representation of 
$G(z,y,y),$ that is,
\begin{equation}
\ln(G(z,y,y))=d(y)+ \int_\bbR d\lambda\,
((\lambda-z)^{-1}-\lambda(\lambda^2
+1)^{-1})\xi(\lambda,y), \lb{4.114}
\end{equation}
to rewrite the trace formula \eqref{4.112} in the form  
 originally obtained in \cite{GHSZ95}, \cite{GS96}, 
\begin{equation}
V(y)=E_0+ \lim_{z\downarrow -\infty} \int_{E_0}^\infty 
d\lambda \, z^2(\lambda-z)^{-2}(1-2\xi(\lambda,y)), 
\quad E_0=\inf (\spec(H_{y,\pi/2})). \lb{4.115}
\end{equation}
Here we used the elementary facts 
\begin{equation}
\xi(\lambda,y)=1-\eta_{y,\alpha_F}(\lambda), \quad 
\xi^{(0)}(\lambda)=1-\eta_{\alpha_F^{(0)}}^{(0)}
(\lambda)=\begin{cases} 1/2, & \lambda >0, \\ 
0, & \lambda < 0. \end{cases} \lb{4.116}
\end{equation}

Of course this formalism is not restricted to the case 
where $\text{def}\,(\dot H)=(1,1).$ The analogous
construction in the case $\text{def}\,(\dot H)=(n,n),$ 
$n\in\bbN,$ then yields an abstract approach to 
matrix-valued trace formulas. This formalism is applicable 
to matrix Schr\"odinger operators and reproduces the 
matrix-valued trace formula analog of \eqref{4.115} 
first derived in \cite{GH97}. To actually prove a 
formula of the type 
\eqref{4.115} for a matrix-valued
$n\times n$ potential $V(x)$ using this abstract 
framework, one then factors
the imaginary part of the 
$n\times n$ matrix $\Im(-G(i,y,y)^{-1})$ into
\begin{equation}
\Im(-G(i,y,y)^{-1})=S(y)^*S(y) \lb{4.1118}
\end{equation}
for some $n\times n$ matrix $S(y)$ and uses relations of 
the type
\begin{equation}
G(z,y,y)=-(S(y)^*M_{y,\alpha_F}(z)S(y) + 
\Re(-G(i,y,y)^{-1}))^{-1}, \quad z\in\bbC\backslash\bbR,  
\lb{4.119}
\end{equation}
where $M_{y,\alpha_F}(z)$ denotes the $n\times n$ 
Donoghue $M$-matrix (cf.~\cite{GKMT98}, \cite{GMT98}, 
\cite{GT97})
for the coresponding matrix-valued  Friedrichs extension
$H_{\alpha_F}$ of $\dot H$ and matrix-valued exponential 
Herglotz representations of the type
\begin{align}
&\ln(-G(z,y,y)^{-1})=C(y)+\int_\bbR d\lambda\,((\lambda 
-z)^{-1}-\lambda(\lambda^2+1)^{-1})\Upsilon(\lambda,y), 
 \lb{4.120} \\
& C(y)=C(y)^*, \quad 0\leq \Upsilon(\lambda,y)\leq I_n 
\text{ for~a.e. } \lambda\in\bbR, 
\quad z\in\bbC\backslash\bbR. \no
\end{align}
Since the actual details are a bit involved, we 
will return to this topic elsewhere \cite{GM99}.

\vspace*{3mm}
\noindent {\bf Acknowledgments.}
F.~G. would like to thank A.~G.~Ramm and P.~N.~Shivakumar 
for their kind invitation  to take part
in this conference, and all organizers for providing a
most stimulating atmosphere during this meeting.


\end{document}